\pdfoutput=1
\documentclass[3p]{elspreprint}

\usepackage{amsthm}
\usepackage[intlimits]{amsmath}
\usepackage{amssymb} 
\usepackage{amsthm}
\usepackage{amsfonts}
\usepackage{bm}
\usepackage{mathdots} %% for \iddots

\usepackage[utf8]{inputenc} % for UTF-8 encoding
\usepackage[english]{babel}  
\usepackage[T1]{fontenc}
\usepackage{times}
\usepackage{mathptmx}  % times roman, including math
\DeclareSymbolFont{cmletters}{OML}{cmr}{m}{n}
\DeclareMathAlphabet{\mathcal}{OMS}{cmsy}{m}{n} % recover 'mathcal'

\journal{Computer Methods in Applied Mechanics and Engineering}
\ifdefined\copyReadyOn
\usepackage{ecrc}
\volume{00}
\firstpage{1}
\runauth{G. Beer, B. Marussig, J. Zechner, Ch. Duenser, T.P. Fries}
\journalname{\journal}
\jid{CMAME}
\jnltitlelogo{CMAME}
\CopyrightLine{2015}{Published by Elsevier Ltd.}
\fi

\usepackage{setspace}
\usepackage{comment}
\usepackage{booktabs}

\providecommand\url[1]{\emph{#1}}

\DeclareMathOperator{\dx}{d}

\newcommand\mal{\,} %% quickly insert \cdot here or leave it as it is
\renewcommand\vec[1]{\mathbf{\boldsymbol{#1}}} %% bold vector
\newcommand\mat[1]{\mathbf{\boldsymbol{#1}}} %% matrix
\newcommand\pt[1]{\mathbf{\boldsymbol{#1}}} %% spatial poiint
\newcommand\ofpt[1]{(\pt{#1})} 
\newcommand\boundary{S} %% Notation for boundary
\newcommand\domain{V} %% notation for omega

\newcommand\fromto[2]{\{ #1, \dots, #2 \}}
\newcommand\range[2]{[#1,#2]}

\newcommand{\figref}[1]{Figure~\ref{#1}}

\newtheoremstyle{myremark}% ⟨name⟩ 
{3pt}% Space above
{3pt}% Space below
{}% Body font
{}% Indent amount
{\bfseries}% Theorem head font % {\itshape}% 
{:}% Punctuation after theorem head
{.5em}% Space after theorem head
{}% Theorem head spec (can be left empty, meaning ‘normal’)
\theoremstyle{myremark}
\newtheorem{remark}{Remark}

\usepackage{listings}
\usepackage{color}
\definecolor{dkgreen}{rgb}{0,0.6,0}
\definecolor{gray}{rgb}{0.5,0.5,0.5}
\definecolor{mauve}{rgb}{0.58,0,0.82}

\begin{document}

%%%%%%%%%%%%%%%%%%%%%%%%%%%%%%%%%%%%%%%%%%%%%%%%%%%%%%%%%%%%%%%%%%%%%%
%%% Title
%%%%%%%%%%%%%%%%%%%%%%%%%%%%%%%%%%%%%%%%%%%%%%%%%%%%%%%%%%%%%%%%%%%%%%

\title{Isogeometric Boundary Element Analysis with elasto-plastic inclusions. Part 1: Plane problems}
\begin{frontmatter}

\author[ifbaddr,newcastleaddr]{Gernot Beer\corref{cor1}}
\author[ifbaddr]{Benjamin Marussig}
\author[ifbaddr]{Jürgen Zechner}
\author[ifbaddr]{Christian Dünser}
\author[ifbaddr]{Thomas-Peter Fries}

\address[ifbaddr]{Institute of Structural Analysis, Graz University
  of Technology, Lessingstraße 25/II, 8010 Graz, Austria}

\address[newcastleaddr]{Centre for Geotechnical and Materials Modelling, University of Newcastle,
  Callaghan, NSW 2308, Australia}

\cortext[cor1]{Corresponding author. 
  Tel.: +43 316 873 6181, fax: +43 316 873 6185, mail: \url{gernot.beer@tugraz.at}, web: \url{www.ifb.tugraz.at}}

\begin{abstract}
In this work a novel approach is presented for the isogeometric Boundary Element analysis of domains that contain inclusions with different elastic properties than the ones used for computing the fundamental solutions. In addition the inclusion may exhibit inelastic material behavior. In this paper only plane stress/strain problems are considered.

In our approach the geometry of the inclusion is described using NURBS basis functions. The advantage over currently used methods is that no discretization into cells is required in order to evaluate the arising volume integrals.  
The other difference to current approaches is that Kernels of lower singularity are used in the domain term. 
The implementation is verified on simple finite and infinite domain examples with various boundary conditions. Finally a  practical application in geomechanics is presented.

\end{abstract}


\end{frontmatter}

% \linenumbers

\section{Introduction}
Isogeometric analysis \cite{Hughes2005a} has gained significant popularity in the last decade because of the fact that geometry data can be taken directly from Computer Aided Design (CAD) programs, potentially eliminating the need for mesh generation. A true companion to CAD is the Boundary Element Method (BEM) because both employ a surface definition of the problem to be solved.

\begin{figure}
\begin{center}
\includegraphics[scale=1.0]{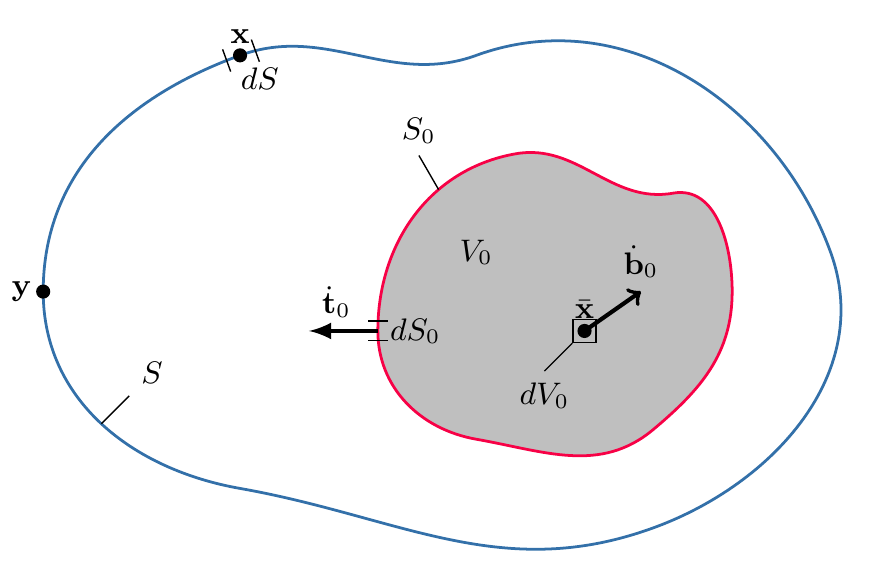}
\caption{Explanation of the derivation of  the integral equation with volume effects}
\label{fig:Betti}
\end{center}
\end{figure}

However, with a  pure surface discretization the BEM can only analyze homogeneous, elastic domains. The method will be extended here to include heterogeneous, inelastic domains by introducing volume effects.
We explain this on an elastic domain with an inclusion $V_{0}$ where body forces are present.
 Using the theorem of Betti as explained in \cite{BeerSmithDuenser2008b}, the boundary integral equation  can be written in incremental form and in matrix notation as:
 \begin{equation}
 \begin{aligned}
   \label{}
   \vec{c} \mal \dot{\vec{u}}\ofpt{y} = \int_{\boundary} \mat{U} (
   \pt{y},\pt{x} ) \dot{\vec{t}}\ofpt{x}
   \dx \boundary  +  \int_{\boundary_{0}} \mat{U}( \pt{y}, \bar{\pt{x}} ) \dot{\vec{t}}_{0}( \bar{\pt{x}} ) \dx \boundary_{0} \\
   - \int_{S} \mat{T}( \pt{y}, \pt{x}) \dot{\vec{u}}\ofpt{x} \dx
   \boundary + \int_{\domain_{0}} \mat{U}( \pt{y},\bar{\pt{x}} )
   \dot{\vec{b}}_{0}\ofpt{x} \dx \domain_{0}
 \end{aligned}
\end{equation}
where \textbf{c} is a free term, $\mat{U}( \pt{y},\pt{x} )$ and $\mat{T}( \pt{y},\pt{x} )$ are matrices containing fundamental solutions for the displacements and tractions at a point $\pt{x}$ due to a source at a point $\pt{y}$ \cite{Beer2015}, $\dot{\vec{u}}\ofpt{x}$ and  $\dot{\vec{t}}\ofpt{x}$ are increments of the displacement and traction vectors on the surface $\boundary$, defining the problem domain (see \figref{fig:Betti}). $\dot{\vec{b}}_{0}( \bar{\pt{x}} )$ are increments of body force inside the inclusion and $\dot{\vec{t}}_{0}( \bar{\pt{x}} )$ are increments of tractions related to the body force acting on surface $\boundary_{0}$ bounding $\domain_{0}$.

\begin{remark}
  In all previous work on elasto-plasticity, integral equations are
  used that involve a higher singularity Kernel in the volume integral
  and the direct use of initial stresses instead of body forces. Here
  we use a different approach involving body forces and a lower
  singularity Kernel. The derivation of the integral equations is
  shown in \ref{Appendix_A}.
\end{remark}

The integral equations can be solved for the unknowns $\vec{u}$ or $\vec{t}$ by discretization. As in majority of previous work on the isogeometric BEM \cite{Simpson2012a, Scott2013197, Marussig2014a, Beer2014a, Marussig2015a, Beer2015a, Beer2015b} we use the collocation method, i.e. we write the integral equations for a finite number $N$ source points $\pt{ y}_{n}$
\begin{equation}
   \label{eq:bie}
   \begin{aligned}
     \vec{c} \mal \dot{\mat{u}}\left( \pt{y}_{n} \right) &=
     \int_{\boundary} \mat{U}\left( \pt{y}_{n},\pt{x} \right) \dot{\vec{t}}\left( \pt{x} \right) \dx \boundary  +  \int_{\boundary_{0}} \mat{U}\left( \pt{y}_{n},\bar{\pt{x}} \right) \dot{\vec{t}}_{0}\left( \bar{\pt{x}} \right) \dx \boundary_{0} \\
      &- \int_{\boundary} \mat{T}\left( \pt{y}_{n},
       \pt{x}\right) \dot{\vec{u}}\left( \pt{x} \right) \dx \boundary
     + \int_{\domain_{0}} \mat{U}\left( \pt{y}_{n},\bar{\pt{x}}
     \right) \dot{\vec{b}}_{0}\left( \bar{\pt{x}} \right) \dx
     \domain_{0}
   \end{aligned}
 \end{equation}
with $n=\fromto{1}{N}$. For the discretization of the surface integrals over $\boundary$ we divide the boundary into patches and use a geometry independent field approximation approach for each patch, i.e. we use different basis functions for the description of the geometry and for the field values.

\begin{equation}
\begin{aligned}
  \pt{x}^{e}  &=  \sum_{k=1}^{K}     R_{k}(u) \mal \pt{x}_{k}^{e} \\
  \vec{u}^{e}  &=  \sum_{k=1}^{K^{d}} R_{k}^{d} (u) \mal \vec{u}_{k}^{e} \\
  \vec{t}^{e}  &=  \sum_{k=1}^{K^{t}} R_{k}^{t} (u) \mal \vec{t}_{k}^{e}
\end{aligned}
\end{equation}
In above equations the superscript $e$ refers to the number of the patch, $R_{k}$, $R_{k}^{d}$ and $R_{k}^{t}$ are NURBS basis functions  with respect to the local coordinate $u$ of the parametrization for the geometry, displacements and tractions respectively. The parameters $\pt{x}_{k}^{e}$ specify the location of control points. The parameters $\vec{u}_{k}^{e}$ and $\vec{t}_{k}^{e}$ are the displacements and tractions mapped to control points. $K$, $K^{d}$, $K^{t}$ specify the number of parameters for each patch.

For an external Neumann problem for example the system of equations
\begin{equation}
\label{SystemOfEquations}
\left[\mathbf{ T}\right]\{\mathbf{ u}\}= \{\mathbf{ F}\} + \{\mathbf{ F}\}_{0}
\end{equation}
is obtained where $\left[\mat{ T}\right]$ is an assembled matrix with coefficients related to Kernel $\mat{T}$ and $\{\vec{ u}\}$ is a vector that collects all displacement components on points $\pt{ y}_{n}$. On the right hand side of \eqref{SystemOfEquations} the vector $\{\vec{ F}\}$ is related to the the given tractions and $\{\vec{ F}\}_{0}= \{\vec{ F}\}^{\boundary_{0}}_{0} + \{\vec{ F}\}^{\domain_{0}}_{0}$ is related to the body force effects, i.e. to the integrals over $\boundary_{0}$ and $\domain_{0}$ in \eqref{eq:bie}.

Details of the implementation of the isogeometric BEM for elastic homogeneous domains can be found in \cite{Beer2015, BeerBordas2014}. Here we concentrate on the definition of the geometry of the inclusions and on the evaluation of the volume integrals.

\section{Basic approach and previous work}

The basic approach is to solve the problem in an iterative way. First the elastic problem is solved considering an elastic homogeneous domain. Then the solution is modified to account for the presence of inclusions and inelastic behavior.

The procedure can be summarized as follows:
\begin{enumerate}
  \item Solve the elastic, homogeneous problem and determine the increment of stress $\dot{\vec{\sigma}}$ inside the inclusion $\domain_0$.
  \item Determine an increment in initial stress $\dot{\vec{\sigma}}_{0}$ due to the fact that the elastic material properties of the inclusion are different from the ones used for the fundamental solutions and/or due to the fact that the elastic limit has been exceeded.
  \item Convert $\dot{\vec{\sigma}_{0}}$  to body force and traction increments $\dot{\vec{ b}}_{0}$, $\dot{\vec{ t}}_{0}$.
  \item Compute new right hand side by evaluating the arising volume and surface integrals.
  \item Solve for the new right hand side and compute a new increment of stress $\dot{\vec{\sigma}}$ inside the inclusion.
  \item Repeat 2. to  5. until $\dot{\vec{\sigma}}_{0}$ is sufficiently small.
\end{enumerate}

\subsection{Elastic inclusions}
\label{ElasI}
Elastic inclusions can be modeled with the multi-region method (see for example \cite{BeerSmithDuenser2008b}) and this involves an additional discretization and increases the number of unknowns. Here we include their treatment in the iterative process required for plasticity.

To compute the initial stress increment for the case where the inclusions have elastic properties which are different to the ones used for the fundamental solutions we use the relation between increments of stress $\dot{\vec{\sigma}}$ and strain $\dot{\vec{\epsilon}}$ in Voigt notation
\begin{eqnarray}
\label{Const}
\dot{\vec{\sigma}} & = & \mat{ C} \mal \dot{\vec{\epsilon}} \\
\dot{\vec{\epsilon}} & = & \mat{ C}^{-1} \mal \dot{\vec{ \sigma}} 
\end{eqnarray}
where \textbf{C} is the constitutive matrix for the domain used for the computation of the fundamental solutions.
The difference in stress between the inclusion and the domain and therefore the initial stress increment can be computed by
\begin{eqnarray}
\dot{\vec{\sigma_{0}}} & = & (\mat{ C}_{i} - \mat{ C}) \mal  \dot{\vec{\epsilon}}
\end{eqnarray}
where $\mat{ C}_{i}$ is the constitutive matrix for the inclusion.

\subsection{Inelastic behavior}
\label{PlasI}
If the inclusion experiences inelastic behavior then additional initial stresses are generated.
Here we use the concept of visco-plasticity, but it is obvious that the method presented here can also be applied to elasto-plasticity.
In visco-plasticity we specify a visco-plastic strain rate
\begin{equation}
\frac{\partial \vec{\epsilon}^{vp}}{\partial t}=\frac{1}{\eta}\Phi(F)\frac{\partial Q}{\partial \vec{ \sigma}}
\end{equation}
where $\eta$ is a viscosity parameter, $F$ is the yield function, $Q$ the plastic potential \cite{simo1998}. It holds that
\begin{eqnarray}
\Phi(F)  =  0 & for & F<0\\
\Phi(F)  =  F & for & F>0.
\end{eqnarray}

The visco-plastic strain increment during a time increment $\Delta t$  can be computed by an explicit scheme
\begin{equation}
\dot{\vec{\epsilon}}^{vp}= \frac{\partial \vec{\epsilon}^{vp}}{\partial t} \mal \Delta t .
\end{equation}
The time step $\Delta t$ can not be chosen freely and if chosen too large, oscillatory behavior will occur in the solution. Suitable time step values can be found in \cite{Cormeau}. The initial stress increment is given by
\begin{equation}
\dot{\vec{\sigma}}_{0}= \mathbf{ C} \mal \dot{\vec{\epsilon}}^{vp}.
\end{equation}

\subsection{Previous work}

The first BEM formulation for inelastic problems has been proposed in \cite{swedlow1971}. The method has been improved substantially in \cite{telles1979}, \cite{telles1991} and \cite{bonnet1996}. The last approach proposed an initial strain formulation in which the consistent tangential operator \cite{simo1998} is used to obtain convergence of quadratic order for the iterative solution procedure.

%% Cells
The common approach for the evaluation of the necessary domain integration is to use cells which are identical to isoparametric finite elements. The cell based method for solving inelastic problems with the BEM is explained in detail in \cite{GaoDavies2002b} and \cite{BeerSmithDuenser2008b}.
To overcome the need for a volume discretization, approaches such as the dual reciprocity BEM \cite{henry1988} or the use of radial basis functions \cite{gao2002} have been proposed. A comparison of these methods to the cell based approach found in \cite{ingber2001} recommends the latter for accuracy and robustness. Moreover, radial basis functions are not suitable for the analysis of infinite domains.
The possibility to automatically generate cells in the inelastic region has been explored in \cite{RibeiroBeerDuenser2008a}. 
In \cite{RiedererDuenserBeer2009a} and \cite{Riederer2010phd} the cell method is extended to cover the simulation of elastic inclusions with different material properties and in \cite{zechner2013} applied to a fast BEM formulation.

All these approaches require the generation of a mesh of cells, which adds an additional effort. Therefore the main innovation presented here is that the concept of cells is abandoned and replaced by a geometry definition of the inclusion as will be explained. It is expected that this approach will not only make the simulation of these problems more user friendly but we expect also an increase in accuracy of the results because the approximation of initial stress inside cells with basis functions is avoided. An important aspect regarding accuracy is that most published cell based methods a continuous variation of initial stresses is assumed, regardless of the fact that the elasto-plastic boundary may cut through a cell and that in this case the variation of the initial stress is discontinuous.

In the following it is first outlined how the geometry of inclusions is defined using NURBS curves and how the arising volume and surface integrals are numerically evaluated.

\section{Geometry definition for inclusions}
\label{Map}

The first task is the description of the geometry of the subdomain $\domain_{0}$. 
For this we propose to use a mapping method introduced recently for trimmed surfaces in \cite{Beer2015a} and \cite{Beer2015}. This means that the subdomain is defined by two NURBS curves and a linear interpolation between them.

\begin{figure}
\begin{center}
\includegraphics[scale=0.5]{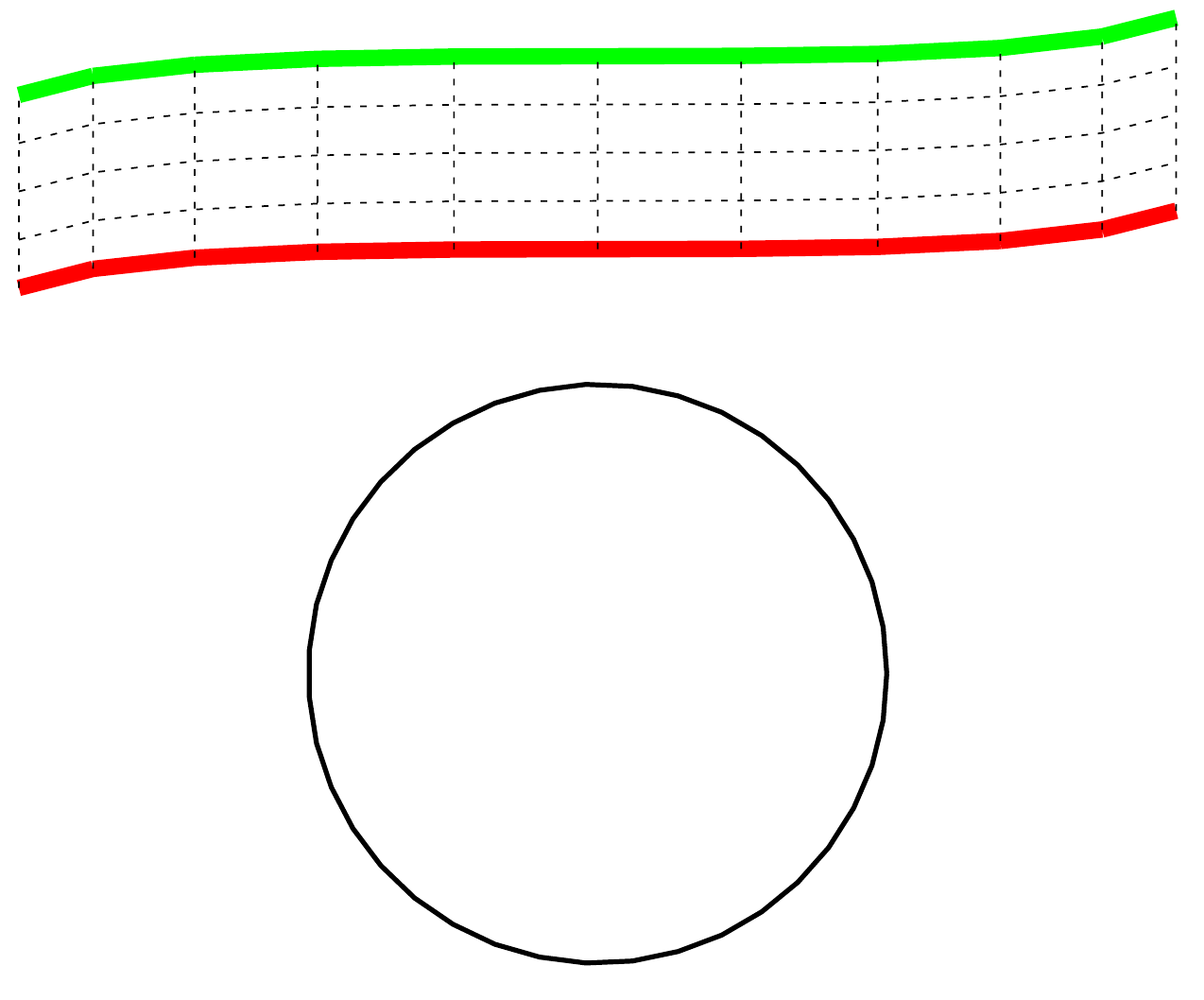}
\includegraphics[scale=0.4]{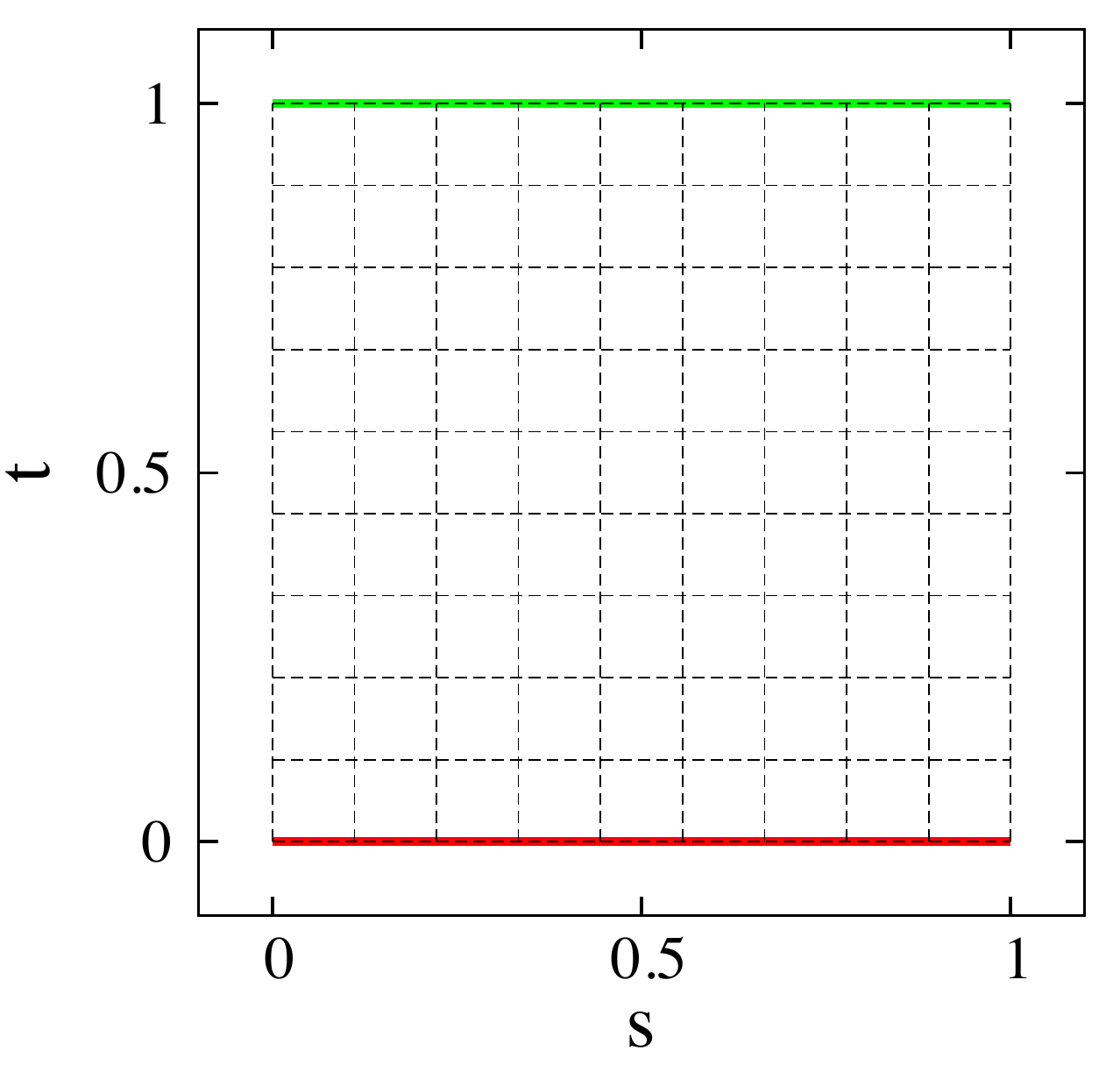}
\caption{A circular excavation with an inclusion above it: Definition of the geometry of an inclusion with 2 NURBS curves in left global and right local coordinate space}
\label{Incl}
\end{center}
\end{figure}

We establish a local coordinate system $\vec{s}=(s,t)^{\intercal}=\range{0}{1}^2$ as shown in \figref{Incl} and perform all computations such as integration and differentiation in this system and then map it to the global $x,y$-system. Note that there is a one to one mapping between the coordinate $s=\range{0}{1}$ and coordinate $u$ of the red and green NURBS curve in \figref{Incl} . The global coordinates of a point $\vec{x}$ with the local coordinates $\pt{s}$ are given by
\begin{equation}
\pt{ x}({s,t})= (1-t) \mal \pt{ x}^{I}(s) + {t} \mal \pt{ x}^{II}({s}) 
\end{equation}
where 
\begin{eqnarray}
  \pt{ x}^{I}({s})=\sum_{k=1}^{K^{I}} R_{k}^{I}({s}) \mal \mathbf{ x}_{k}^{I} & \text{and} &\pt{ x}^{II}({s})=\sum_{k=1}^{K^{II}} R_{k}^{II}({s}) \mal \pt{ x}_{k}^{II} .
\end{eqnarray}
The superscript $I$ relates to the bottom (red) curve and $II$ to the top (green) curve and $ \mathbf{ x}_{k}^{I} $, $ \mathbf{ x}_{k}^{II} $ are control point coordinates. $K^{I}$ and $K^{II}$ are the number of control points, $R_{k}^{I}({s})$ and $R_{k}^{II}({s})$ are NURBS basis functions.
The derivatives are given by
\begin{equation}
\begin{aligned}
  \frac{\partial \pt{ x}({s,t})}{\partial {s}}&=& (1-{t}) \mal \frac{\partial \pt{ x}^{I}({s})}{\partial {s}} + {t} \mal \frac{\partial \pt{ x}^{II}({s}) }{\partial {s}} \\
  \frac{\partial \pt{ x}({s,t})}{\partial {t}}&=& -\pt{ x}^{I}({s}) +
  \mal \pt{ x}^{II}({s})
\end{aligned}
\end{equation}
where

\begin{equation}
\begin{aligned}
  \frac{\partial\pt{ x}^{I}({s})}{\partial {s}}&=&\sum_{k=1}^{K^{I}} \frac{\partial R_{k}^{I}({s})}{\partial {s}} \mal \pt{ x}_{k}^{I} \\
  \frac{\partial\pt{ x}^{II}({s})}{\partial {}s}&=&\sum_{k=1}^{K^{II}}
  \frac{\partial R_{k}^{II}({s})}{\partial {s}} \mal \pt{ x}_{k}^{II}
  .
\end{aligned}
\end{equation}
The Jacobian matrix of this mapping is 
\begin{equation}
\label{jacobian1}
\mat{J}=
\begin{pmatrix}
  \frac{\partial x}{\partial {s}} & \frac{\partial y}{\partial {s}}  \\ \\
  \frac{\partial x}{\partial {t}}  & \frac{\partial y}{\partial {t}}
\end{pmatrix}
\end{equation}
and the Jacobian is $J\ofpt{s}=| \mathbf{ J} |$.

\begin{remark}
  It is obvious that the following theory is not restricted in any way
  to the simple geometry description outlined above. Any method that
  allows the mapping of the geometry to a unit square can be
  applied. In the following we have used this description
  mainly in order to simplify the explanation of the method.
\end{remark}

\section{Computation of $\{\vec{ F}\}_{0}$}
Here we discuss the computation of the right hand side during iteration. This involves the evaluation of the integrals in Equation (\ref{eq:bie}) over $\boundary_{0}$ and $\domain_{0}$ using Gauss quadrature, similar in a way to the Nyström method proposed in \cite{Zechner2015}.

\subsection{Computation of the surface integral over $\boundary_{0}$}
\label{S0}
In order to minimize the number of Gauss points a subdivision into integration regions is recommended. Such subdivision is essential if the plastic zone does not extend over the whole domain, because in this case the integrand will be discontinuous.

For the computation of the integral over surface $\boundary_{0}$ we first integrate along the two curves defining the inclusion (i.e. for $u=\range{0}{1}$) and then along the edges (i.e. from $t=\range{0}{1}$).
For the integration along the bounding curves the global locations of Gauss points are computed by
\begin{eqnarray}
\pt{ x}^{I}({\mathrm{u}})=\sum_{k=1}^{K^{I}} R_{k}^{I}({u}) \mal \pt{ x}_{k}^{I} & \text{and} &\pt{ x}^{II}({u})=\sum_{k=1}^{K^{II}} R_{k}^{II}({u}) \mal \pt{ x}_{k}^{II} 
\end{eqnarray}
The Jacobian of this transformation is $J^{i}= \sqrt{(\frac{dx^{i}}{du})^2 + (\frac{dy^{i}}{du})^2}$.

Gauss integration requires the use of a local coordinate system $\xi=\range{-1}{1}$.
The transformation from $u$ to $\xi$ coordinates is given by
\begin{equation}
u  = \frac{\Delta u_{n_s}}{2} (1+\xi) + u_{n_s} 
\end{equation}
where $\Delta u_{n_s}$ is the size of the integration region $n_s$ and $u_{n_s}$ is the start coordinate.
The Jacobian is $J_{u}^{n_s}=\frac{\Delta u_{n_s}}{2}$.

For the integration along the left edge $e_1$ we have
\begin{equation}
  \pt{ \bar{x}}({t}) = (1-{t}) \mal \pt{ x}^{I}_{1} + {t} \mal \pt{ x}^{II}_{1} 
\end{equation}
where $\pt{ x}^{I}_{1}$ and $\pt{ x}^{II}_{1}$ are the coordinates of a point  on the top and bottom curves with the local coordinate $u=0$. 

Assuming for simplicity\footnote{This will be the case for the examples presented later which involve thin inclusions, but is not a restriction of the method.} that there is no subdivision in the $t$ direction the transformation to $\xi$ coordinates is given by
\begin{equation}
t= \frac{1}{2}(1+\xi).
\end{equation}
The Jacobian of this transformation is $J_{e_1}=\frac{1}{2}(\pt{ x}^{II}_{1} - \pt{ x}^{I}_{1})$ 

For the right edge $e_2$ we have
\begin{equation}
\pt{ \bar{x}}({t}) = (1-{t}) \mal \pt{ x}^{I}_{2} + {t} \mal \pt{ x}^{II}_{2} 
\end{equation}
where $\pt{ x}^{I}_{2}$ and $\pt{ x}^{II}_{2}$ are the coordinates of a point on the top and bottom curves respectively for $u=1$. The Jacobian of this transformation is $J_{e_2}=\frac{1}{2}(\pt{ x}^{II}_{2} - \pt{ x}^{I}_{2})$. 
We can now write the sub vector of $\{\mathbf{ F}\}_{0}^{S_{0}}$ related to the collocation point $n$ as
\begin{equation}
  \label{GaussS}
   \mathbf{ F}^{S_{0}}_{0n}=    \sum_{i=1}^2 \sum_{n_s=1}^{N^s}\int_{-1}^{1}\mat{U}( \pt{y}_{n},\bar{\pt{x}} ) \dot{\pt{t}}_{0}( \bar{\pt{x}} ) \mal J^{i} \mal J_{u}^{n_s} \mal  \mal \dx \xi + \sum_{j=1}^2  \int_{-1}^{1}\mat{U}( \pt{y}_{n},\bar{\pt{x}} ) \dot{\pt{t}}_{0}( \bar{\pt{x}} )  \mal J_{e_j} \mal  \mal \dx \xi .
\end{equation}
where $N^s$ is the number of subregions.
Applying Gauss integration we have for example
\begin{equation}
  \label{Gaussint}
   \sum_{i=1}^2 \sum_{m=1}^{M} \pt{U}\left( \pt{y}_{n},\bar{\pt{x}}(\xi_{m}) \right) \dot{\pt{t}}_{0}\left( \bar{\pt{x}}(\xi_{m}) \right) \mal J_{e_i}  \mal W_{m} 
\end{equation}
for the second term in \eqref{GaussS} where $M$ is the number of integration points, $\xi_{m}$ are the local coordinates of Gauss points and $W_{m} $ are quadrature weights. To determine the number of Gauss points necessary for an accurate integration we consider that whereas there is usually a moderate variation of body force the Kernel $\mat{U}$ behaves like $\mathcal{O}(\ln r)$ with $r=| \pt{y}_n-\bar{\pt{x}} |$ and approaches infinity as $\bar{\pt{ x}}$ approaches $\pt{ y}_{n}$, so the number of integration points has to be increased if $\boundary_{0}$ is close to the boundary $\boundary$. 

If the collocation point is located on $\boundary_{0}$ the integrand approaches infinity as the point is approached. Here we apply the method that is used for dealing with weakly singular integrals over surface $S$. It involves the transformation to a local coordinate system $\gamma=\range{-1}{1}$ where the Jacobian tends to zero as the collocation point is approached.
Details of the implementation can be found in \cite{Beer2015}.

\subsection{Computation of the volume integral over $V_{0}$}

For the volume integration the transformation from $\vec s$ coordinates to  $\vec{\xi}=(\xi,\eta)^{\intercal}=\range{-1}{1}^2$ is  given for integration region $n_s$ by
\begin{eqnarray}
s & = \frac{\Delta s_{n_s}}{2} (1+\xi) + s_{n_s} \\
\nonumber
t & = \frac{\Delta t_{n_s}}{2} (1+\eta) + t_{n_s}
\end{eqnarray} 
where $\Delta s_{n_s}\times \Delta t_{n_s}$ denotes the size of the integration region and $s_{n_s},t_{n_s}$ are the starting coordinates. The Jacobian of this transformation is for the integration over $\boundary_{0}$ is {$J_{\xi}^{n_s}=\frac{\Delta s_{n_s} \mal  \Delta t_{n_s}}{4}$}.

The sub vector of $\{\vec{ F}\}_{0}^{\boundary_{0}}$ related to collocation point n can be written as:
\begin{equation}
  \label{Gauss}
   \vec{ F}^{V_{0}}_{0n} = \sum_{n_s=1}^{N^s}\int_{-1}^{1}  \int_{-1}^{1} \mat{U} \left( \pt{y}_{n},\bar{\pt{x}}(\xi,\eta) \right)
  \dot{\vec{b}}_{0} \left( \bar{\pt{x}} (\xi,\eta) \right) J\ofpt{s} \mal J_{\xi}^{n_s} \mal \dx \xi \dx \eta
\end{equation}

Applying Gauss integration we have:
\begin{equation}
  \label{GaussintV}
   \mathbf{ F}^{V_{0}}_{0n}  \approx \sum_{n_s=1}^{N^s}\sum_{m=1}^{M} \sum_{n=1}^{N} \mat{U}\left( \pt{y}_{n},\bar{\pt{x}}(\xi_{m},\eta_{n}) \right) \dot{\vec{b}}_{0}\left( \bar{\pt{x}}(\xi_{m},\eta_{n}) \right) J\ofpt{s}  \mal J_{\xi}^{n^s} \mal W_{m} \mal W_{n}
\end{equation}
where $N^s$ is the number of integration regions and $M$ and $N$ are the number of integration points in $\xi$ and $\eta$ directions respectively. To determine the number of Gauss points necessary for an accurate integration we consider that whereas there is usually a moderate variation of body force, the Kernel \textbf{U} is $O(\ln r)$ so the number of integration points has to be increased if $\pt{y}_n$ is close to $V_0$.

If the integration region includes the collocation point $\textbf{y}_{n}$, then the integrand tends to infinity as the point is approached and a procedure has to be invoked that has been used to deal with weakly singular integrals in three-dimensional BEM, involving triangular subregions.

\begin{figure}
  \begin{center}
    \includegraphics[scale=0.55]{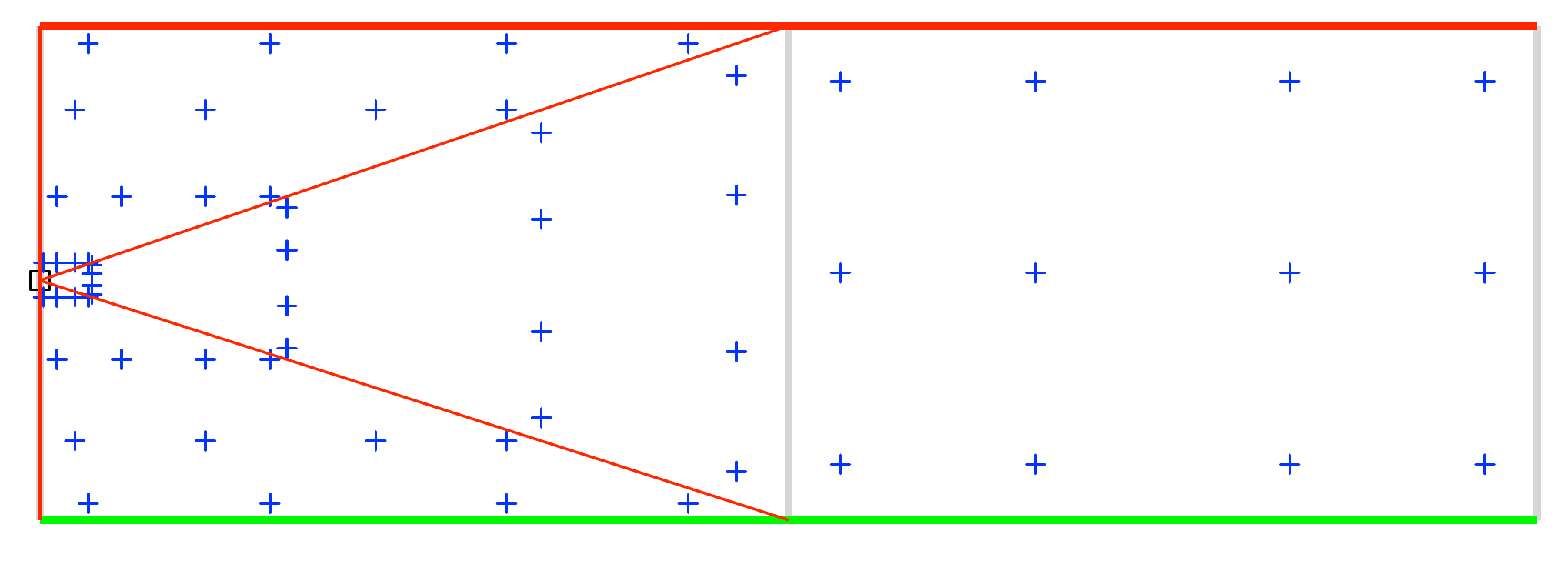}
    \caption{Subdivision into integration regions when the collocation point (marked by a square) is part of $\domain_{0}$. The grey line indicates a subdivision into 2 integration regions due to the extent of the plastic zone, red thin lines the subdivision into triangular subregions. Location of Gauss points are marked with crosses}
    \label{Sin}
  \end{center}
\end{figure}

In this approach we perform the integration in a local coordinate system, where the Jacobian tends to zero as the singularity point is approached. For this we divide the integration region into two, three or four triangular sub-regions depending on whether the collocation point is at a corner, edge or inside. The procedure is well documented in \cite{Beer2015} and leads to the following expression:
  \begin{equation}
  \label{GaussintVS}
  \begin{aligned}
    \mathbf{ F}^{V_{0}}_{0n}  &=  \sum _{n_s=1}^{N^s} \sum_{k=1}^{K^\triangle} \int_{-1}^{+1} \int_{-1}^{+1} \mat{U}\left( \pt{ y}_{n},\pt{ x} \right) \dot{\vec{b}}_{0}\left( \bar{\pt{x}}(\xi,\eta) \right)  \mal J\ofpt{s}  \mal J_{\triangle_k} \mal J_\xi^{n_s} \mal  \dx \xi \dx \eta\\
    &\approx \sum _{n_s=1}^{N^s} \sum_{k=1} ^{K^\triangle}\sum_{m=1}^{M}
    \sum_{\ell=1}^{L} \mat{U}\left( \pt{ y}_{n},\pt{ x} \right) \:
    \dot{\vec{b}}_{0}\left( \bar{\pt{x}}(\xi_{m},\eta_{\ell}) \right)
    \mal J\ofpt{s} \mal J_{\triangle_k} \mal J_\xi^{n_s} \mal W_{m}
    W_{l}
  \end{aligned}
\end{equation}
where $K^\triangle$ is the number of triangles and $J_{\triangle_k} $ is the Jacobian of the transformation from the square to the triangular subregion. An example of this subdivision is shown in \figref{Sin}.

\subsection{Computation of $\dot{\vec{ t}}_0$ and
  $\dot{\vec{ b}}_0$}

For the evaluation of the integrals \eqref{Gaussint}, \eqref{GaussintV} and \eqref{GaussintVS}, the values of $\dot{\vec{ t}}_0$ and $\dot{\vec{ b}}_0$ must be known at the Gauss points and their determination involves stress or strain evaluations inside the inclusions. Since the location of quadrature points vary depending on the location of the collocation point $\pt{ y}_n$, this would involve a very large number of evaluations for each iteration. In addition some computations would involve singular integration.

It is therefore convenient to compute the required values at regular grid points inside the inclusion and then interpolate or extrapolate the values to the required locations. In the simplest case a linear interpolation can be used. The advantage of this scheme is that differentiations can be carried out numerically using finite differences. For the analysis a regular grid of points is therefore established in the the local $\vec{s}$ coordinate system of the inclusion.

After computing the initial stress increment $\dot{\vec{\sigma}}_{0}$ using the procedures outlined in sections~\ref{ElasI} and \ref{PlasI}, the initial traction increments $\dot{\vec{ t}}_{0}$ are computed from the initial stresses by
% \begin{eqnarray}
%   \dot{t}_{0x}= n_{x}\dot{\sigma}_{0x} + n_{y}\dot{\tau}_{0xy} \\
%   \nonumber
%   \dot{t}_{0y}= n_{y}\dot{\sigma}_{0y} + n_{x}\dot{\tau}_{0xy} 
% \end{eqnarray}
\begin{equation}
  \dot{\vec{t}}_0=
  \begin{pmatrix}
    \dot{\sigma}_{0x} & \dot{\tau}_{0xy} \\
    \dot{\tau}_{0xy}  & \dot{\sigma}_{0y} \\
  \end{pmatrix}
  \vec{n}
\end{equation}
where $\vec{n}$ is the unit outward normal vector to the surface $\boundary_{0}$. The body force increment $\dot{\vec{b}}_{0}$ can be computed by
% \begin{eqnarray}
%   \dot{\mathrm{b}}_{0x} & = & -(\frac{\partial \dot{\sigma}_{0x}}{\partial x} + \frac{\partial \dot{\tau}_{0xy}}{\partial y}) \\
%   \nonumber
%   \dot{\mathrm{b}}_{0y} & = & -(\frac{\partial \dot{\sigma}_{0y}}{\partial y} + \frac{\partial \dot{\tau}_{0xy}}{\partial x}) 
% \end{eqnarray}
\begin{equation}
  \dot{\vec{b}}_0 = -
  \begin{pmatrix}
    \frac{\partial\dot{\sigma}_{0x}}{\partial x} + 
    \frac{\partial\dot{\tau}_{0xy}}{\partial y} \\ \\
    \frac{\partial\dot{\tau}_{0xy}}{\partial x} +
    \frac{\partial\dot{\sigma}_{0x}}{\partial y} 
  \end{pmatrix}
\end{equation}
It is convenient to compute the derivatives with respect to local coordinates \textbf{s} first and then transform them to global coordinates. For example the global derivatives of $\sigma_x$ in terms of local derivatives are given by the transformation
\begin{equation}
  \label{Sigg}
  \vec{\sigma}_{x,\vec{ x}}= \vec{ J}^{-1} \mal \vec{\sigma}_{x,\vec{s}}
\end{equation}
where
\begin{eqnarray}
  \vec{\sigma}_{x,\vec{ x}}  = 
  \begin{pmatrix}
    \frac{\partial \sigma_{x}}{\partial x} \\ \\
    \frac{\partial \sigma_{x}}{\partial y} 
  \end{pmatrix} & \text{and} & \vec{\sigma}_{x,\vec{ s}} =
                               \begin{pmatrix}
                                 \frac{\partial \sigma_{x}}{\partial s} \\ \\
                                 \frac{\partial \sigma_{x}}{\partial t}
                               \end{pmatrix}
\end{eqnarray}
and $\vec{J}$ is the Jacobian matrix \eqref{jacobian1}.

The derivatives are numerically computed using finite differences. For grid points inside the inclusion that have other points left and right (or top and bottom) of them we use a central finite difference, whereas for points that only have one point on a side we use forward or backward finite differences.

\begin{figure}
  \begin{center}
    \includegraphics[scale=0.55]{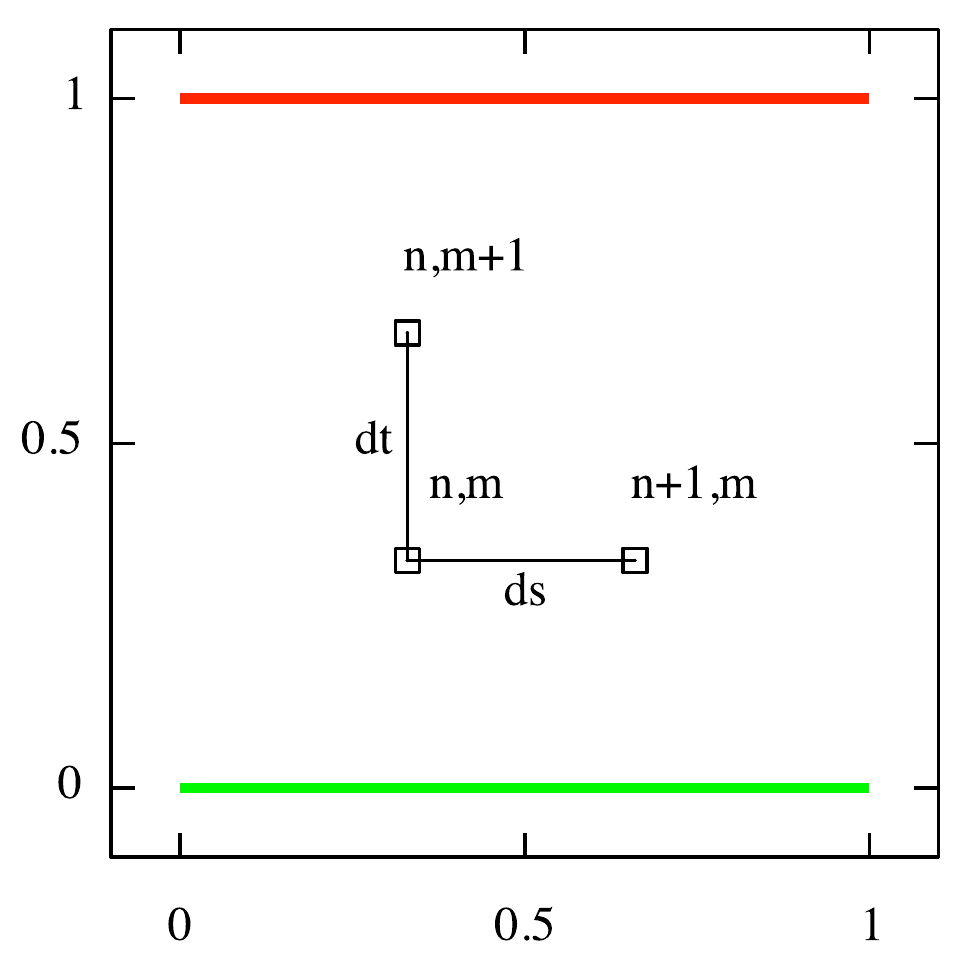}
    \caption{Template for the computation of first derivatives inside
      an inclusion by forward difference in the local s,t coordinate
      system for grid point n,m.}
    \label{FD}
  \end{center}
\end{figure}

Referring to \figref{FD} the local derivatives, computed using a forward finite difference scheme, are given as
\begin{eqnarray}
  \frac{\partial \sigma}{\partial {s}} = \frac{\sigma_{n+1,m} - \sigma_{n,m}}{ds} & \text{and} &
  \frac{\partial \sigma}{\partial {t}} = \frac{\sigma_{n,m+1} - \sigma_{n,m}}{dt} .
\end{eqnarray}

\begin{remark}
  In the following examples we apply a simple linear interpolation between grid points and a linear extrapolation from grid points to the boundary $\boundary_{0}$. Obviously more sophisticated schemes may be applied. However care has to be taken that the variation of the initial stress may be discontinuous. The simple scheme applied here leads to an increase in the accuracy for determination of the derivatives and for the evaluation of the associated integrals as the number of grid points is increased. The convergence of the solution as a function of the number of internal points is investigated in one of the numerical examples below.
\end{remark}

\section{Computation of results inside the inclusion}

As mentioned above, the solution algorithm requires the evaluation of strains and stresses at internal points.  For the initial solution (without body forces) the displacements at a point $\mathbf{ y}_{i}$ inside the inclusion $\domain_0$ can be computed by
\begin{eqnarray}
  \label{eq:int}
  \vec{u}\left( \pt{y}_{i} \right)  = 
  \int_{S} \mat{U}\left( \pt{y}_{i},\pt{x} \right) \vec{t}\left( \pt{x} \right) \dx \boundary  - \int_{S}
  \mat{T}\left( \pt{y}_{i}, \pt{x}\right) \vec{u}\left( \pt{x} \right) \dx \boundary .
\end{eqnarray}
The strain at a point $\mathbf{ y}_{i}$ inside the inclusion is
\begin{eqnarray}
  \vec{\epsilon} \left( \pt{y}_{i} \right)  = 
  \int_{S} \mat{S}\left( \pt{y}_{i},\pt{x} \right) \vec{t}\left( \pt{x} \right) \dx \boundary   -   \int_{S}
  \mat{R}\left( \pt{y}_{i}, \pt{x}\right) \vec{u}\left( \pt{x} \right) \dx \boundary 
\end{eqnarray}
where $\vec{S}$ and $\vec{R}$ are derived fundamental solutions \cite{BeerSmithDuenser2008b}. The fundamental solution $\vec{S}$ has a singularity of order $\mathcal{O}({r}^{-1})$ and $\vec{R}$ a singularity of $\mathcal{O}({r^{-2}})$. The stresses can be computed using \eqref{Const}.

For the subsequent solution (including body forces) the displacements at a point $\pt{ y}_{i}$ inside the inclusion is
\begin{equation}
  \begin{aligned}
    \vec{u}\left( \pt{y}_{i} \right) = \int_{S} \mat{U}\left(
      \pt{y}_{i},\pt{x} \right) \vec{t}\left( \pt{x} \right) \dx
    \boundary - \int_{S}
    \mat{T}\left( \pt{y}_{i}, \pt{x}\right) \vec{u}\left( \pt{x} \right) \dx \boundary\\
    + \int_{S_{0}} \mat{U}\left( \pt{y}_{i},\bar{\pt{x}} \right)
    \dot{\vec{t}}_{0}\left( \bar{\pt{x}} \right) \dx \boundary_{0} +
    \int_{V_{0}} \mat{U}\left( \pt{y}_{i},\bar{\pt{x}} \right)
    \dot{\vec{b}}_{0}\left( \bar{\pt{x}} \right) \dx \domain_{0}
  \end{aligned}
\end{equation}
The strain can be computed by
\begin{equation}
  \label{}
  \begin{aligned}
    \vec{\epsilon} \left( \pt{y}_{i} \right) = \int_{S} \mat{S}\left(
      \pt{y}_{i}, \pt{x} \right) \vec{t}\left( \pt{x} \right) \dx
    \boundary -
    \int_{S} \mat{R}\left( \pt{y}_{i}, \pt{x} \right) \vec{u}\left( \pt{x} \right) \dx \domain \\
    + \int_{S_{0}} \mat{S}\left( \pt{y}_{i},\bar{\pt{x}} \right)
    \dot{\vec{t}}_{0}\left( \bar{\pt{x}} \right) \dx \boundary_{0} +
    \int_{V_{0}} \mat{S}\left( \pt{y}_{i},\bar{\pt{x}} \right)
    \dot{\vec{b}}_{0}\left( \bar{\pt{x}} \right) \dx \domain_{0}
  \end{aligned}
\end{equation}
Again, for the evaluation of the integrals Gauss integration is used. The efficient and accurate evaluation of the integrals over the surface $\boundary$ has been described in some detail in \cite{Beer2015}. Here we concentrate on the evaluation of the integrals over $\boundary_{0}$ and $\domain_{0}$.

\subsection{Computation of integral over $S_{0}$}
Since the result points are inside the inclusion, regular integration as outlined in section \ref{S0} can be used with the number of integration points chosen depending on the proximity of point $\mathbf{ y}_{i}$ to the boundary $S_{0}$.

\subsection{Computation of integral over $V_{0}$}

For the computation of the volume term for the internal point $\pt{ y}_{i}$ we have to check if the point is inside the integration region or not.  In the case it is not, we use normal Gauss integration with the number of integration points depending on the proximity of $\pt{ y}_{i}$ to the integration region. In the case it is, we have to consider that the Kernel $\vec{S}$ tends to infinity as $\mathcal{O}(r^{-1})$ as $\bar{\pt{x}}$ approaches $\pt{y}_{i}$ and is therefore strongly singular.

\begin{figure}
  \begin{center}
    \includegraphics[scale=0.55]{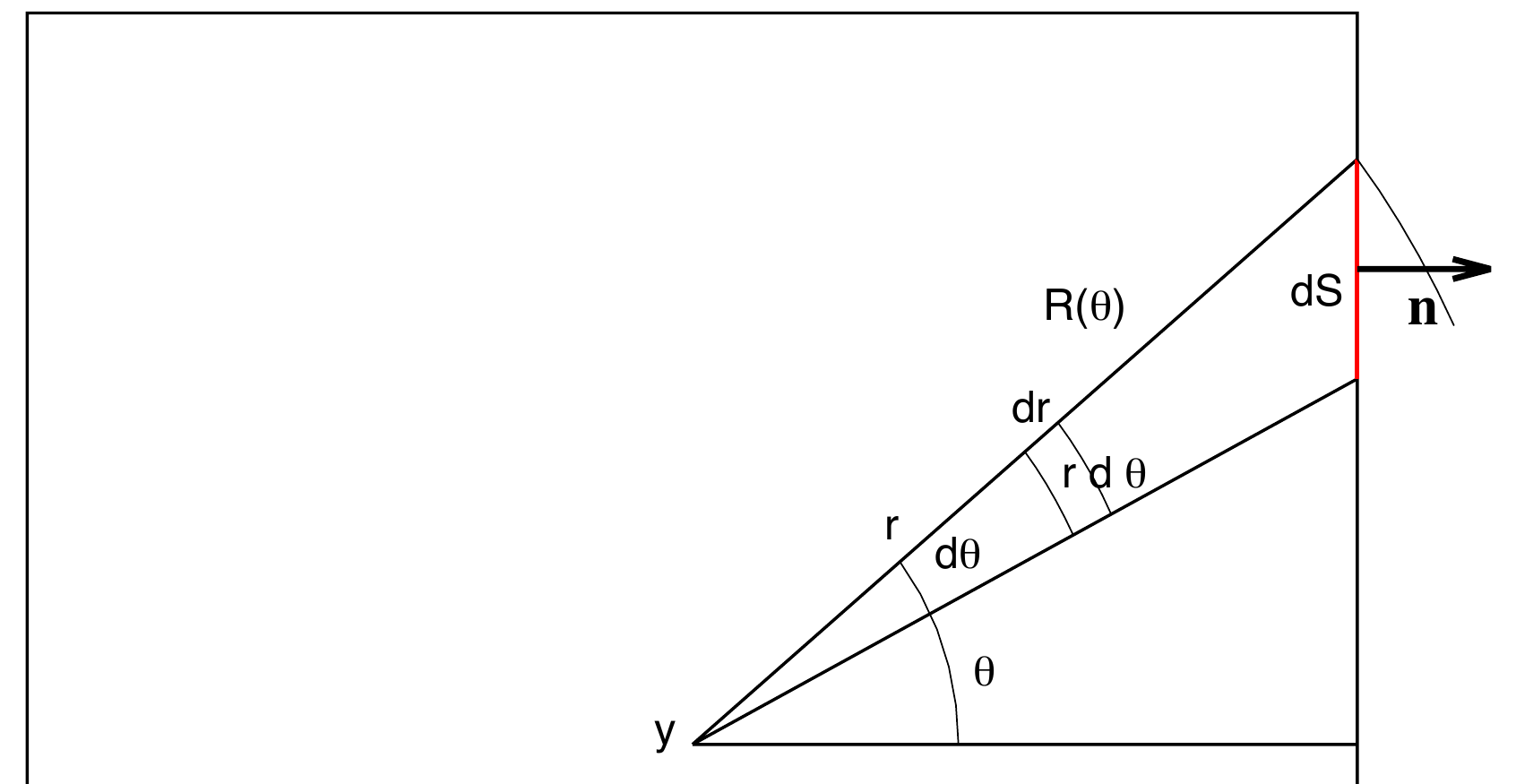}
    \caption{Evaluation of strongly singular integral of S using polar
      coordinates}
    \label{StrSin}
  \end{center}
\end{figure}

The singularity can be isolated by replacing the integral
\begin{equation}
  \begin{aligned}
    \int_{\domain_{0}} \mat{S}\left( \pt{y}_{i},\bar{\pt{x}} \right)
    \dot{\vec{b}}_{0}\left( \bar{\pt{x}} \right) \dx \domain_{0} &=
    \int_{\domain_{0}} \mat{S}\left( \pt{y}_{i},\bar{\pt{x}} \right) \left[ \dot{\vec{b}}_{0}\left( \bar{\pt{x}} \right) - \dot{\vec{b}}_{0}\left( \pt{y}_{i} \right)\right] \dx \domain_{0} \\
    &+ \left[\int_{\domain_{0}} \mat{S}\left( \pt{y}_{i},\bar{\pt{x}}
      \right) dV_{0}\right] \mal \dot{\vec{b}}_{0} \left( \pt{y}_{i}
    \right)
  \end{aligned}
\end{equation}
The first integral is now weakly singular and can be treated as described above. The second integral can be transformed into a surface integral. Using polar coordinates as shown in \figref{StrSin} and defining a small circular area of exclusion around the singularity with a radius of $\epsilon$, we have:
\begin{equation}
\begin{aligned}
  \int_{\domain_{0}} \mat{S}\left( \pt{y}_{i},\bar{\pt{x}} \right) \dx
  \domain_{0} & =
  \int_{0}^{2\pi} \int_{\epsilon}^{R(\theta)}\frac{1}{r} \hat{\mat{ S}} \mal r  \mal \dx \theta  \mal \dx r\\
  & = \int_{0} ^{2\pi} \left[\int_{\epsilon}^{R(\theta)} \dx r \right]
  \hat{\mat{ S} }\mal \dx\theta = \int_{0} ^{2\pi}\hat{ \mat{ S}} \mal
  R(\theta) \mal \dx\theta - \epsilon \int_{0} ^{2\pi}\hat{ \mat{ S}}
  \mal \dx\theta
\end{aligned}
\end{equation}
where $\hat{\mathbf{ S}}$ is the nonsingular part of Kernel $\mat{S}$. It can be shown that the second integral is zero \cite{GaoDavies2002b}.
After substitution of 
\begin{equation}
  R(\theta) \dx\theta = \mathbf{ n} \bullet \mathbf{ r}\frac{1}{R(\theta)} dS
\end{equation}
 we obtain
\begin{equation}
  \label{SdV}
  \int_{\domain_{0}} \mat{ S} \dx \domain_{0} = \int_{\domain_{0}} \frac{1}{R(\theta)} \hat{\mat{ S}} \mal \vec{ n} \bullet
  \mat{ r} \mal \dx \boundary_{0} = \int_{\boundary_{0}} \mat{ S} \mal \vec{ n} \bullet \vec{ r} \mal \dx \boundary_{0}
\end{equation}
where $\boundary_{0}$ is the boundary surface of the inclusion. In equation~\eqref{SdV} $\vec{n}$ is the outward normal and $\vec{r}$ is the position vector.  The volume integral is now replaced by
\begin{equation}
  \int_{\domain_{0}} \mat{S}\left( \pt{y}_{i},\bar{\pt{x}} \right) \left[\dot{\vec{b}}_{0}\left( \bar{\pt{x}} \right) - \dot{\vec{b}}_{0}\left( \pt{y}_{i} \right)\right] \dx \domain_{0} + \left[ \int_{\boundary_{0}} \mat{S} \mal \vec{ n}\bullet \vec{ r} \mal \dx \boundary_{0} \right] \dot{{\vec{b}}}_{0} \left( \pt{y}_{i} \right) .
\end{equation}

\section{Test examples}
In the following sections the theory is tested on simple examples, where the solution is known.

\subsection{Test example 1: Cube with elastic inclusion}
\begin{figure}
  \begin{center}
    \includegraphics[scale=0.6]{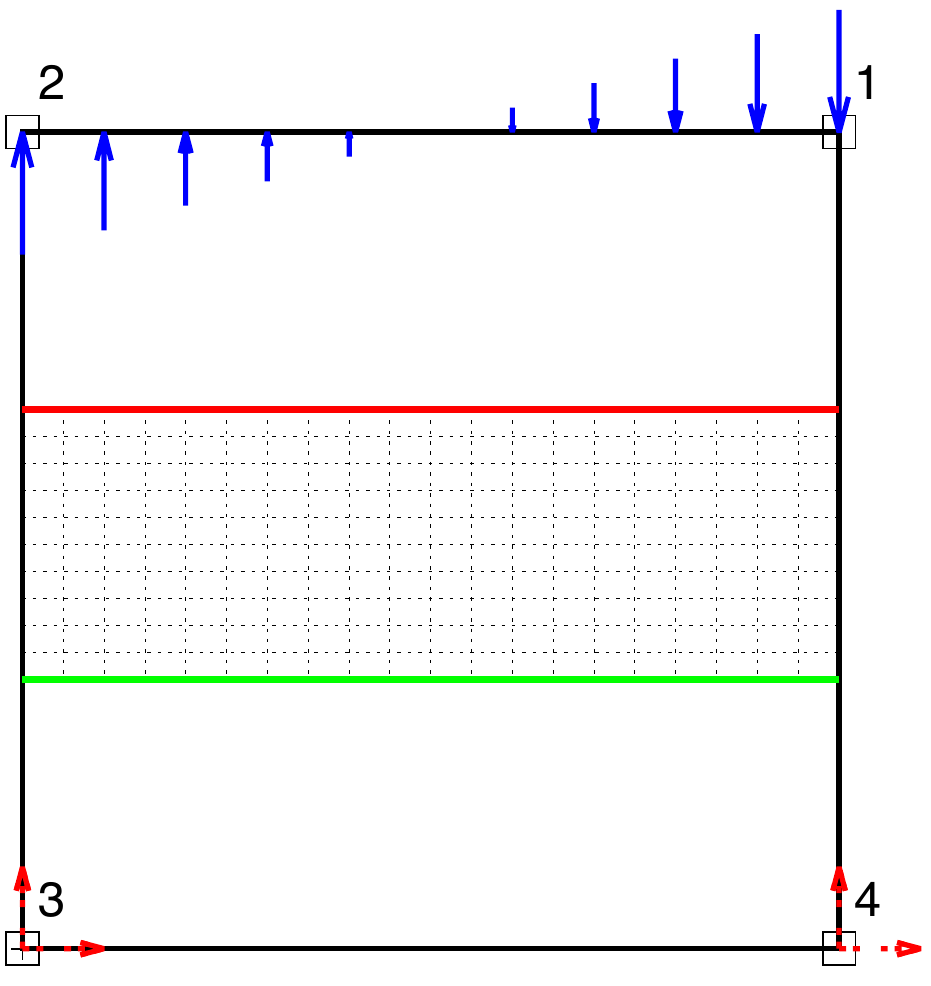}
    \includegraphics[scale=0.6]{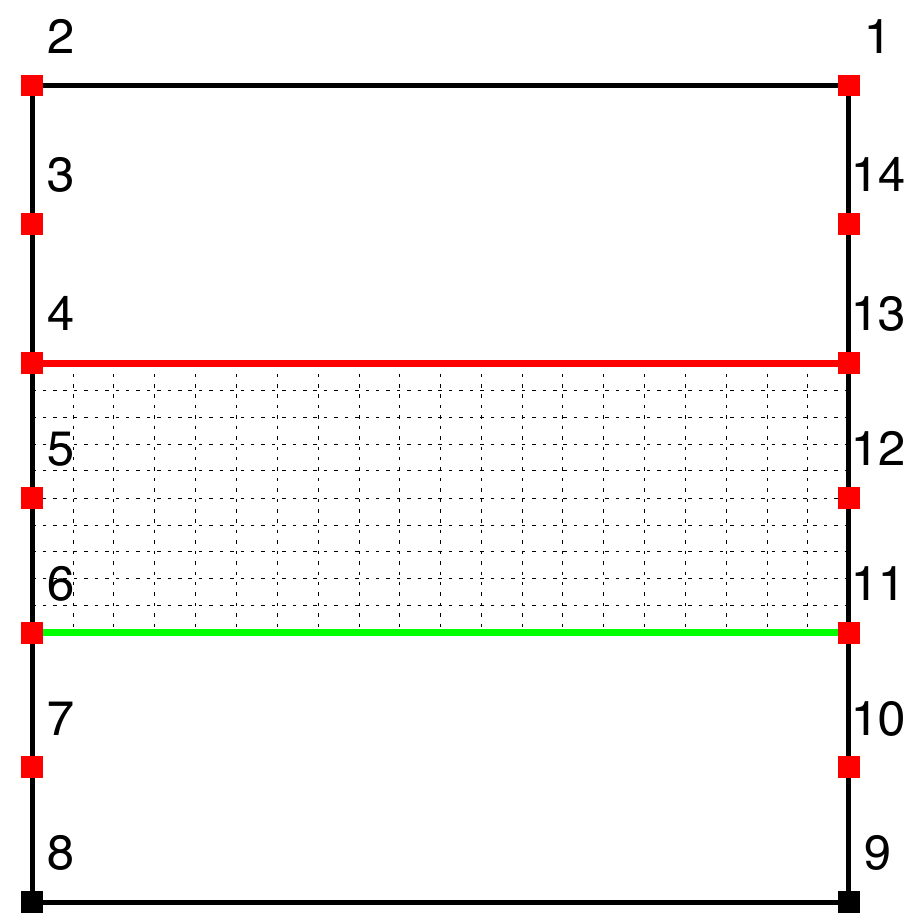}
    \caption{Test example 1. Left: Geometry definition with 4 NURBS
      patches and control points, loading and boundary conditions;
      Right: Collocation points (restrained ones marked in black).}
    \label{Test1}
  \end{center}
\end{figure}

The first example tests the algorithm for the case of a single elastic inclusion. It consists of a cube with the dimension $1\times 1 \times 1$ composed of two different materials. A two-dimensional analysis is carried out using plane stress assumptions and the discretization is shown in \figref{Test1}.
The cube is defined by four NURBS curves with knot vectors and coefficients ($x$, $y$, $z$ coordinates and weights) as follows:
\begin{lstlisting}
Curve 1: Knot vector= 0 0 1 1; Coefficients= 1 1 0 1; 0 1 0 1
Curve 2: Knot vector= 0 0 1 1; Coefficients= 0 1 0 1; 0 0 0 1
Curve 3: Knot vector= 0 0 1 1; Coefficients= 0 0 0 1; 1 0 0 1
Curve 4: Knot vector= 0 0 1 1; Coefficients= 1 0 0 1; 1 1 0 1
\end{lstlisting}
\begin{figure}
  \begin{center}
    \includegraphics[scale=0.47]{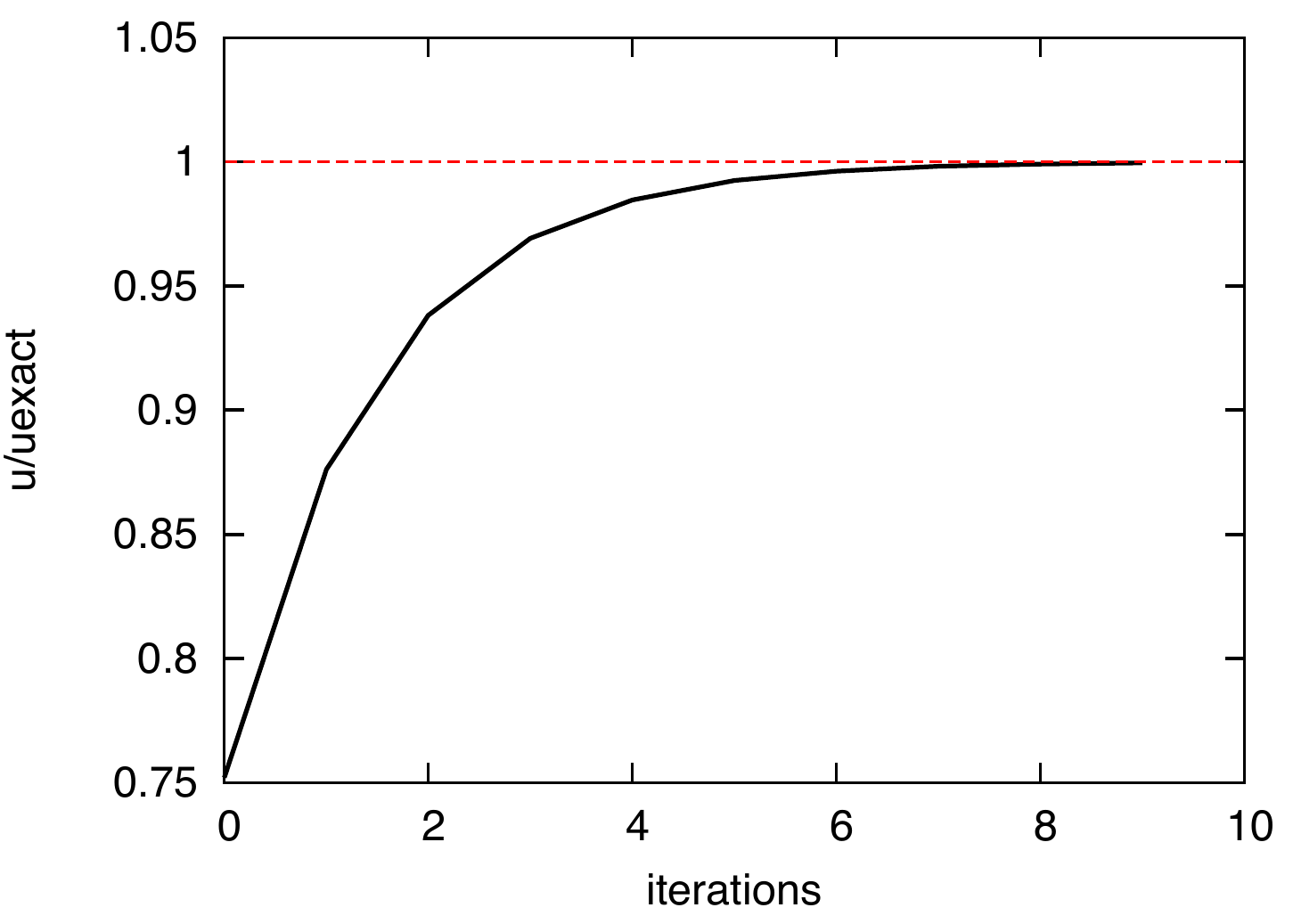}
    \caption{Test example 1: Plot of ratio of maximum computed displacement to the exact one as function of the number of iterations. }
    \label{Test1C}
  \end{center}
\end{figure}
The inclusion is defined by two NURBS curves;
\begin{lstlisting}
Curve 1: Knot vector= 0 0 1 1; Coefs= 1 0.66 0 1; 1 0.66 0 1
Curve 2: Knot vector= 0 0 1 1; Coefs= 0 33 0 1; 1 0.33 0 1
\end{lstlisting}
and assigned a Young's modulus $E$ of half the one used for computing the fundamental solution and no change in the Poisson's ratio $\nu$.  The cube is loaded with a moment as shown left in \figref{Test1} and is fixed at the bottom.

For the analysis, the concept of a geometry independent field approximation was used and the basis functions for approximating the displacements were defined using the following knot vectors
\begin{lstlisting}
Curve 1:  0 0  1 1; 
Curve 2:  0 0 0 0.34 0.34 0.67 0.67 1 1 1; 
Curve 3:  0 0  1 1 ; 
Curve 4:  0 0 0 0.33 0.33 0.66 0.66 1 1 1; 
\end{lstlisting}
with all weights equal to one.
\begin{figure}
  \begin{center}
    \includegraphics[scale=0.4]{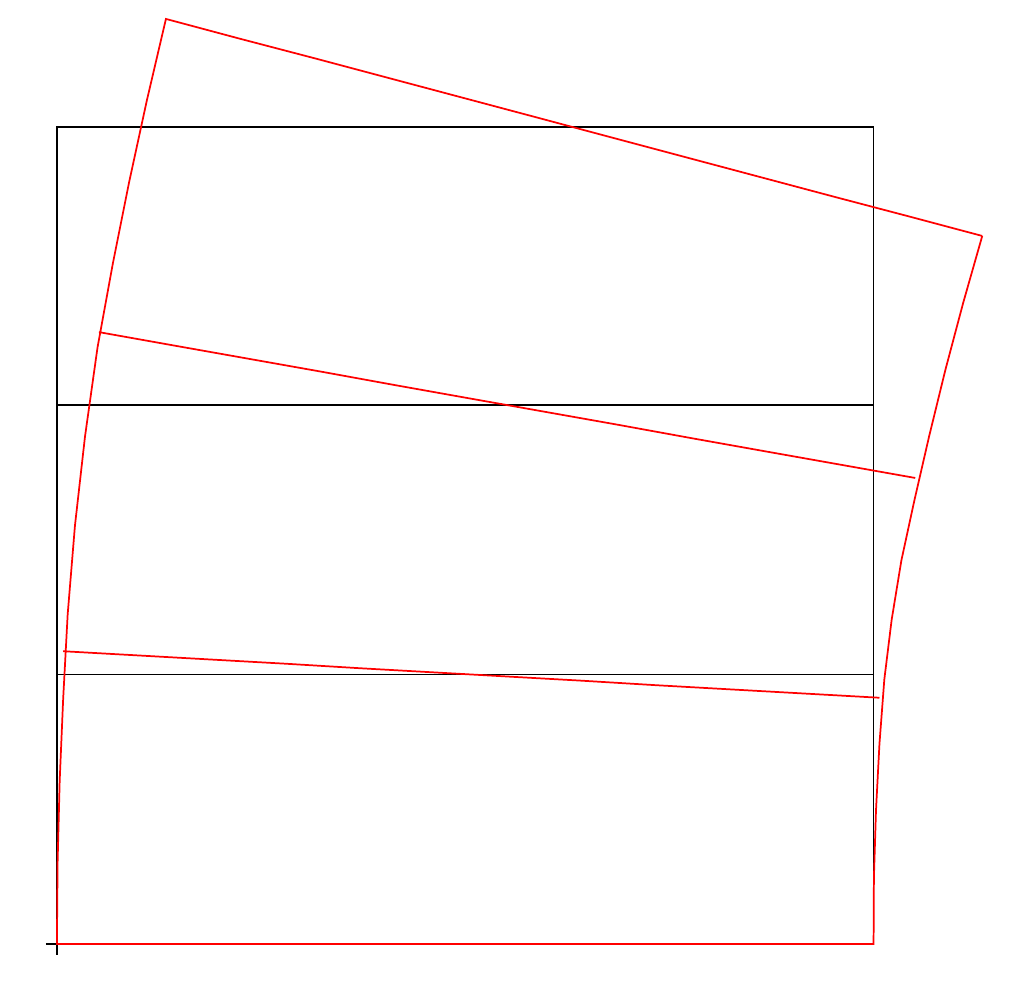}
    \caption{Test example 1: Deformed shape. }
    \label{Test1dis}
  \end{center}
\end{figure}
This approximation results in a quadratic variation of the displacements in the vertical direction with a $C^{1}$ discontinuity at the interface between materials. The resulting location of the collocation points are shown in \figref{Test1} on the right and this will give the exact solution for the applied loading.  The results are shown in \figref{Test1C} and \figref{Test1dis}. It can be seen that convergence to the exact solution is achieved after about $7$ iterations.

\subsection{Test example 2: Cube with inelastic inclusion}

\begin{figure}
  \begin{center}
    \includegraphics[scale=0.45]{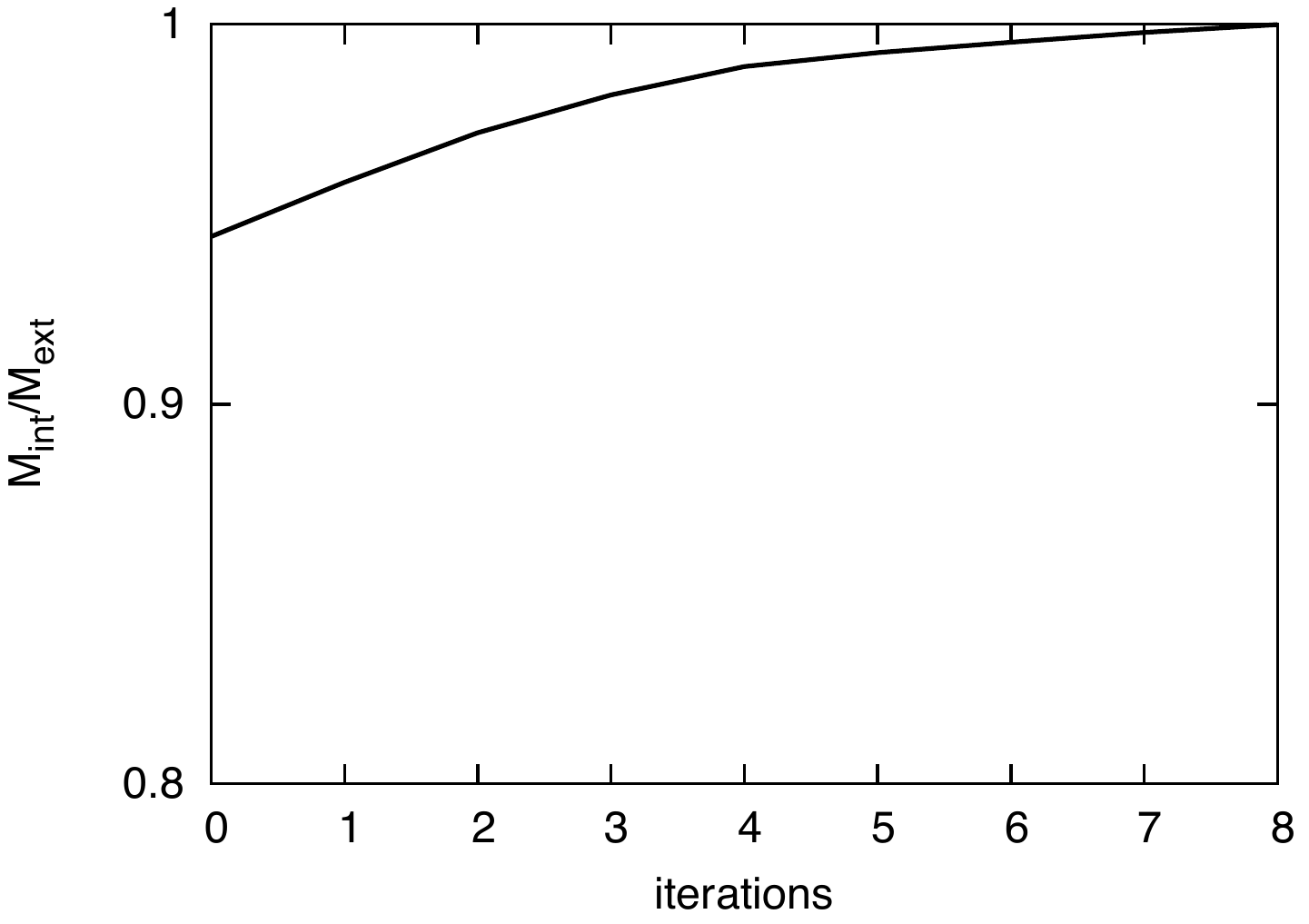}
    \caption{Test example 2: Plot of the ratio internal and external
      moment versus number of iterations}
    \label{conv3}
  \end{center}
\end{figure}

This is the same as example~1, except that this time we assume that the elastic properties are the same everywhere, but inside the inclusion the material behaves in an inelastic way. A \textit{von Mieses} type material law was implemented which restricts the maximum compressive or tensile normal stress to 80\% of the maximum compressive or tensile elastic normal stress.  At the end of the iteration it is checked that the yield condition is satisfied everywhere. An indication of convergence is that the internal moment is equal to the external moment. \figref{conv3} shows the convergence of the iteration based on this criterion after 8 iterations.
\begin{figure}
  \begin{center}
    \includegraphics[scale=0.6]{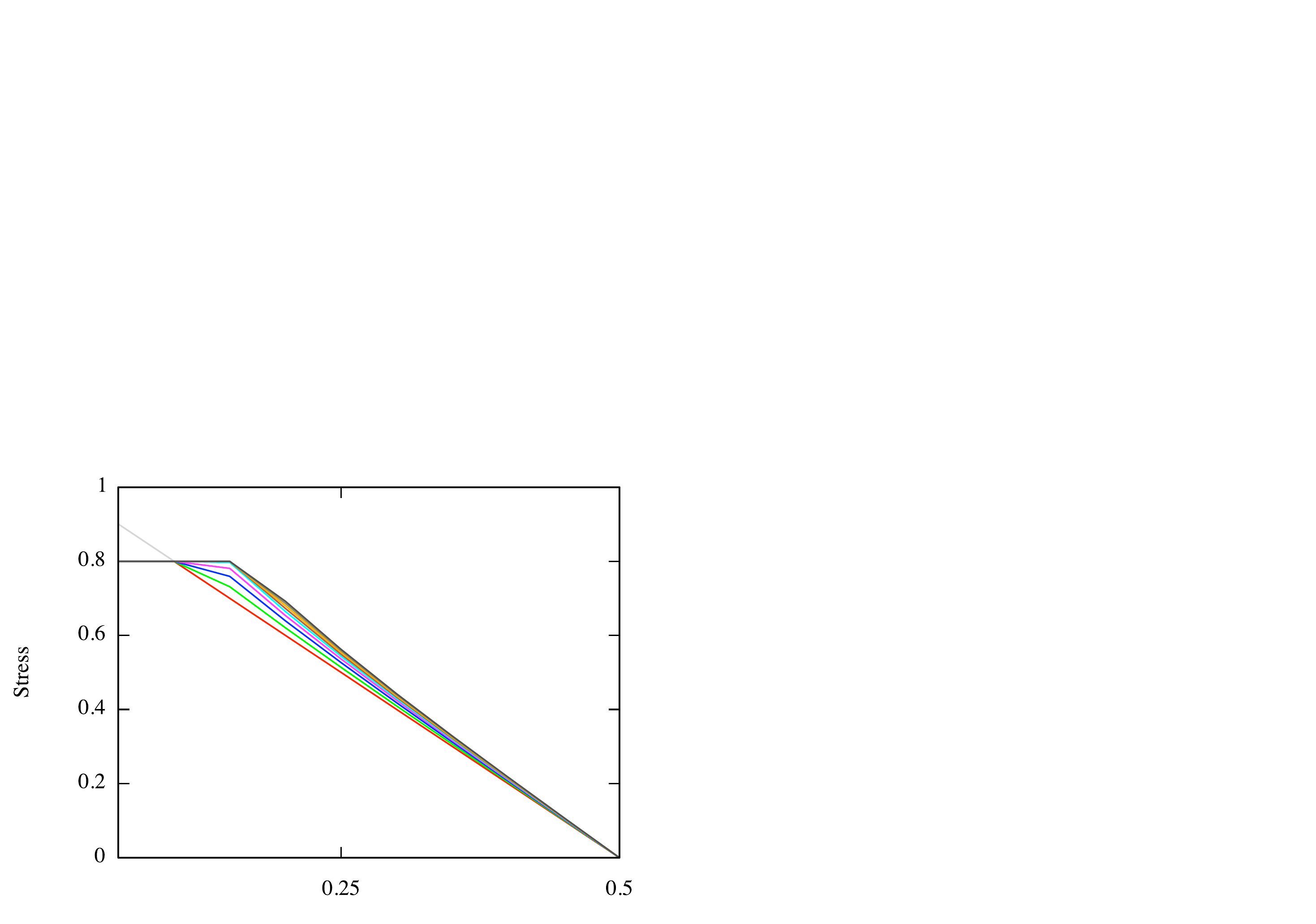}
    \caption{Test example 2: Change in the distribution of normal
      stress on left half of the cube during iterations. Initial
      (elastic) stress state in grey, final stress state in black}
    \label{stress3}
  \end{center}
\end{figure}
\figref{stress3} shows the history of stress changes.

\subsection{Test example 3: Hole in an infinite domain}

For the previous examples $\dot{\mathbf{ b}}_0$ was zero everywhere. Test example 3 is designed to test the volume integration involving this term. It is a hole in an homogeneous infinite domain ($E=100, \nu=0$) under plane strain conditions with an inelastic inclusion on top of it.
\begin{figure}
  \begin{center}
    \includegraphics[scale=0.5]{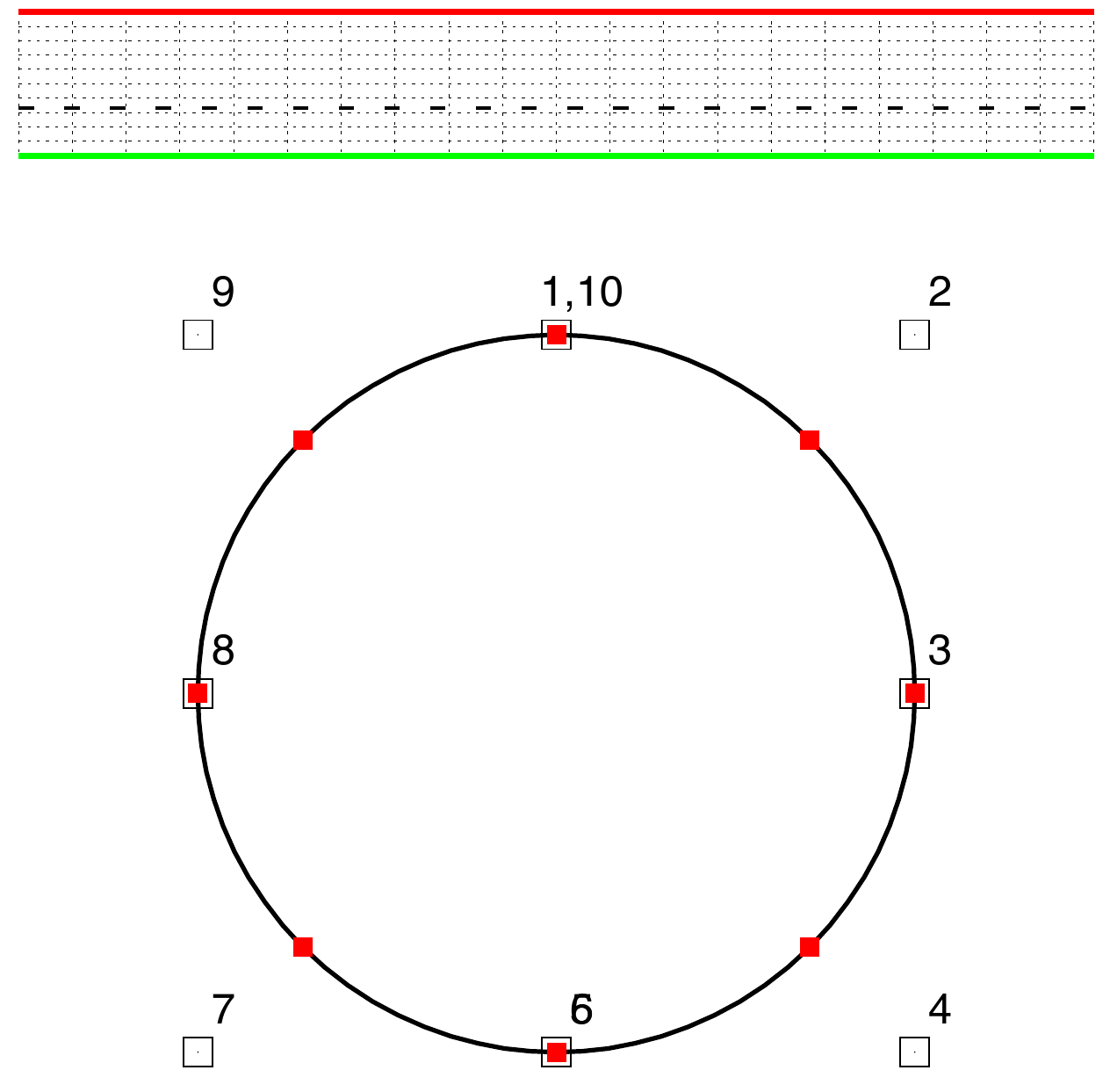}
    \caption{Test example 3: The geometry of the circle is defined with two NURBS patches. Control points are shown as hollow squares, collocation points as filled red squares. The inclusion above the hole is defined by two linear NURBS curves. Also shown is the line (dotted) along which the stresses are plotted.}
    \label{4geo}
  \end{center}
\end{figure}
The geometry of the problem is shown in \figref{4geo}. The hole is defined by two NURBS curves as follows:
\begin{lstlisting}
Curve 1: 
Knot vector= 0 0 0 0.5 0.5 1 1 1;
Coefs= 0.5 1 0 1;1 1 0 0.707;1 0.5 0 1;1 0 0 0.707; 0.5 0 0 1
Curve 2: 
Knot vector= 0 0 0 0.5 0.5 1 1 1
Coefs= 0.5 0 0 1;0 0 0 0.707;0 0.5 0 1;0 1 0 0.707;0.5 1 0 1
\end{lstlisting}
This exactly describes the geometry of a circle. An isogeometric approach was used with the same basis functions defining the approximation of the displacements.

\begin{figure}
  \begin{center}
    \includegraphics[scale=0.6]{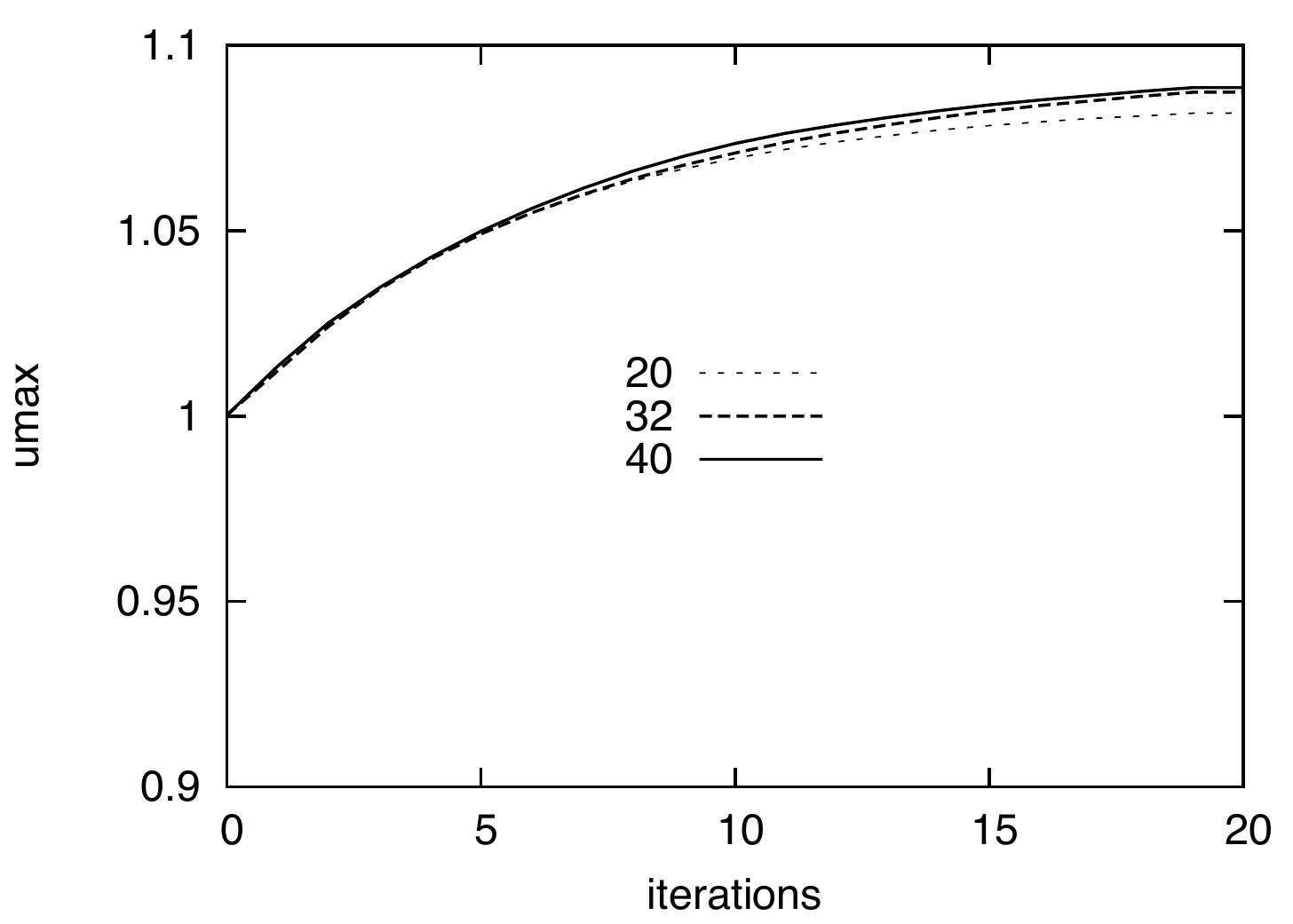}
    \caption{Test example 3: Convergence of top displacement for
      different number of grid points.}
    \label{4conv}
  \end{center}
\end{figure}

\begin{figure}
  \begin{center}
    \includegraphics[scale=0.5]{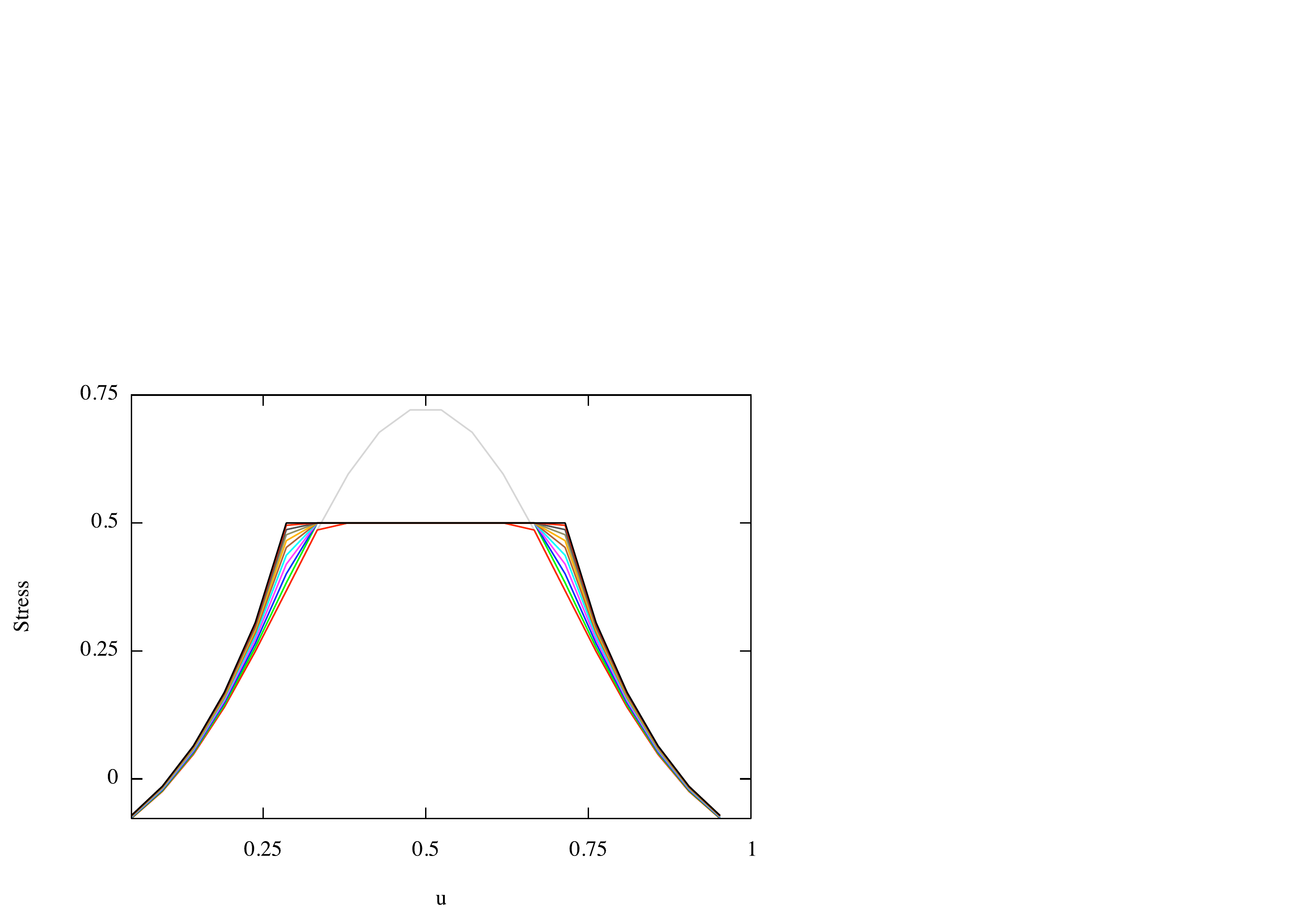}
    \caption{Test example 3: Evolution of vertical stress during
      iteration. Grey line shows the initial elastic stress and black
      line the final stress.}
    \label{4stress}
  \end{center}
\end{figure}

\begin{figure}
  \begin{center}
    \includegraphics[scale=0.4]{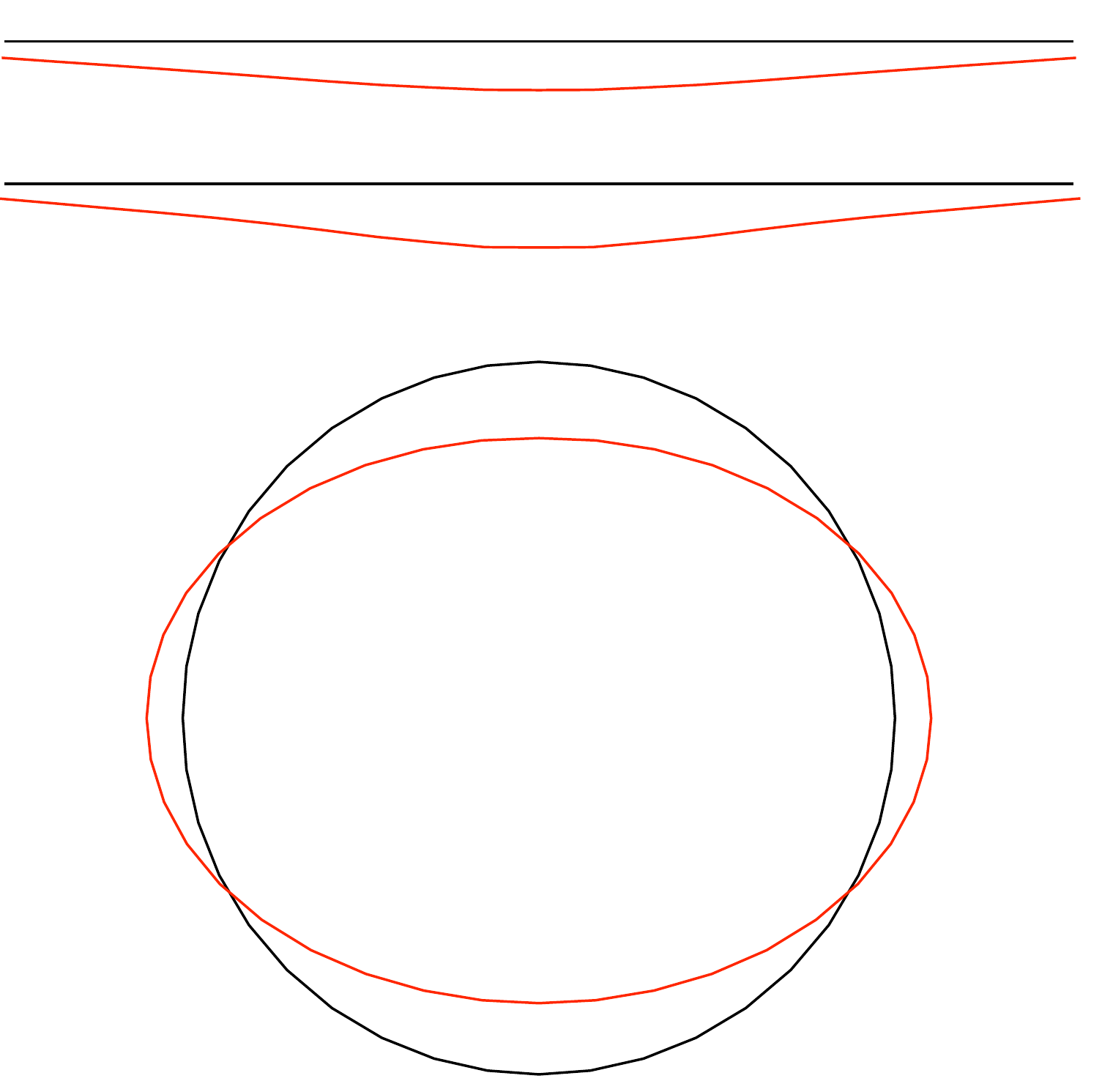}
    \caption{Test example 3: Displaced shape.}
    \label{4disp}
  \end{center}
\end{figure}

The inclusion was defined by two NURBS curves as follows:
\begin{lstlisting}
Curve 1: 
Knot vector= 0 0 1 1 
Coefs= -0.25 1.45 0 1; 1.25 1.45 0 1
Curve 2: 
Knot vector= 0 0 1 1
Coefs= -0.25 1.25 0 1; 1.25 1.25 0 1
\end{lstlisting}

The hole is subjected to a boundary traction along $\boundary$ of $\vec{t}=(0,\,n_{y})^\intercal$, where $n_{y}$ is the vertical component of unit outward normal to $\boundary$. The yield condition limits the tensile stress in the vertical direction to $0.5$. This example was also used to test the sensitivity of the results to the number of grid points. \figref{4conv} shows the convergence of the maximum displacement (at the top of the circle) for the case of $20$ to $40$ points. Results with a higher number of points were indistinguishable from each other.
\figref{4stress} shows the evolution of the vertical stress along the line depicted in \figref{4geo} and \figref{4disp} the displaced shape.

\section{Practical example}

\begin{figure}
  \begin{center}
    \includegraphics[scale=0.55]{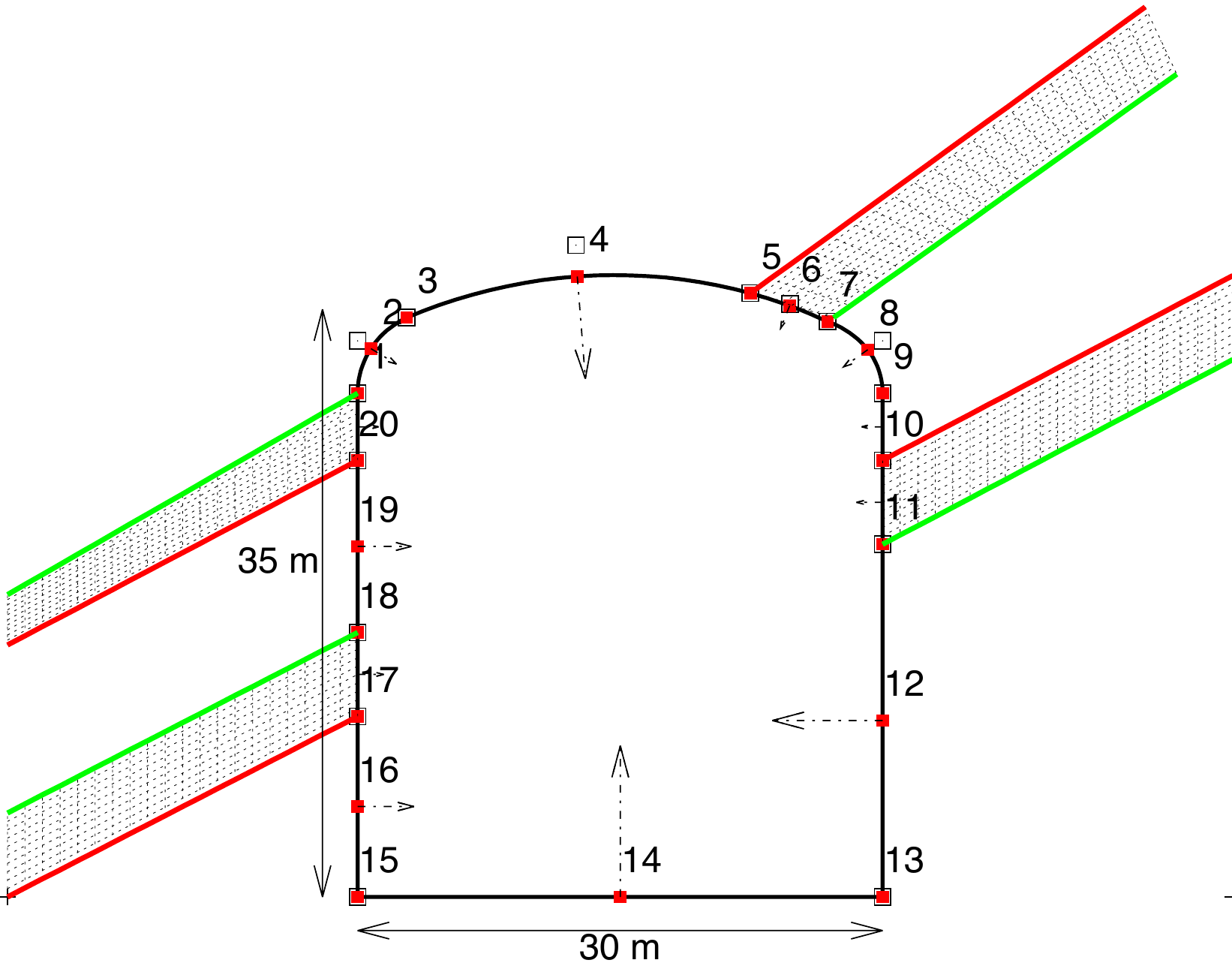}
    \caption{Practical example: Geometry of cavern showing the NURBS patches describing the boundary with control points shown as hollow squares. Four inclusions are described by linear patches. Also shown are the collocation points (red filled squares) used for the analysis.}
    \label{Masjedgeo}
  \end{center}
\end{figure}

The practical example is one that has been solved with a coupled BEM/FEM method and reported in \cite{BeerSmithDuenser2008b}. It relates to the plane strain analysis of an underground power station cavern. \figref{Masjedgeo} shows the geometry of the cavern with the discretization into 12 NURBS patches. The data for the geometry definition are shown in \ref{Appendix_B}.
\begin{figure}
  \begin{center}
    \includegraphics[scale=0.55]{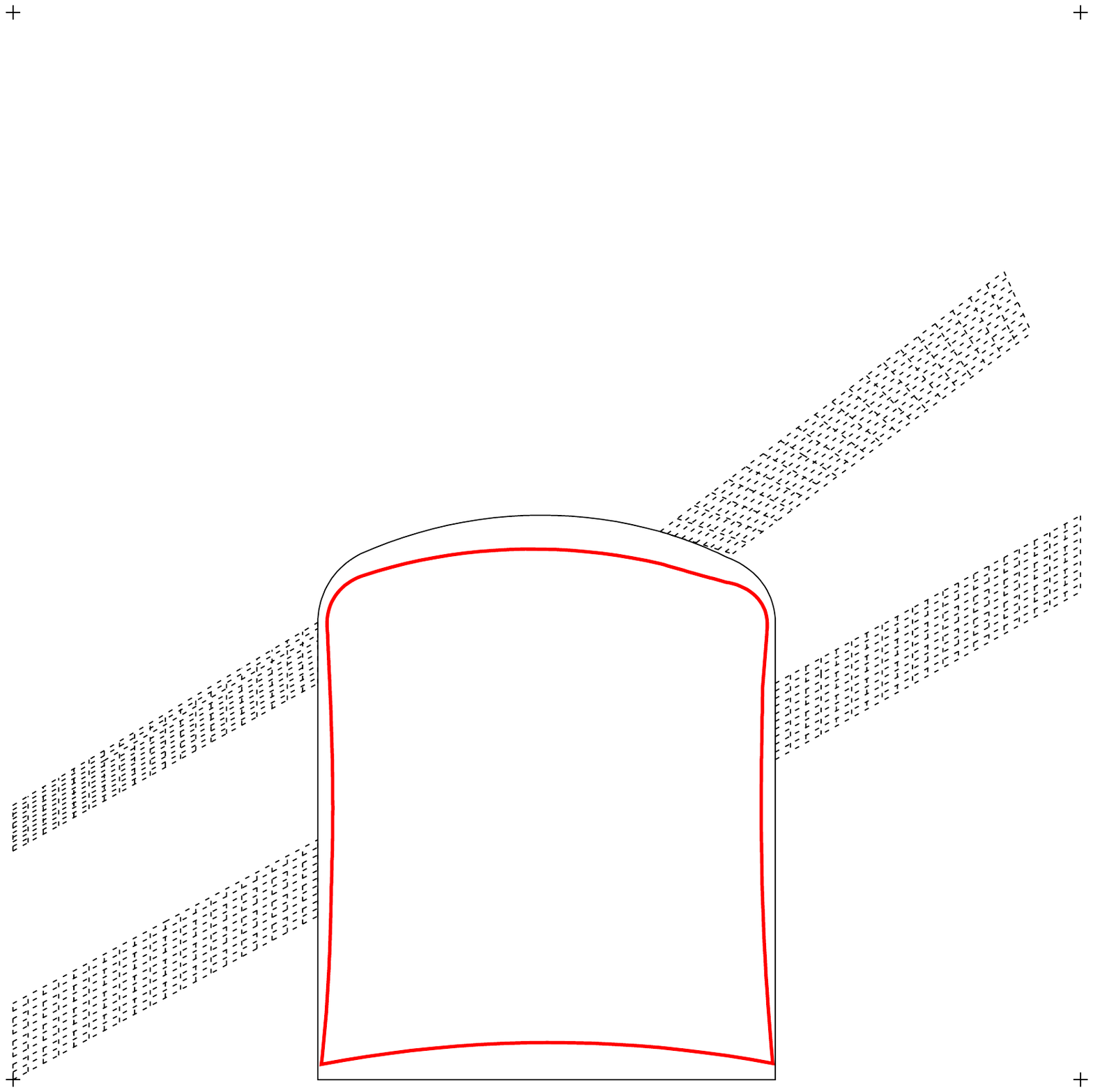}
    \caption{Practical example: Displaced shape.}
    \label{Masjedres}
  \end{center}
\end{figure}

The basis functions used for the description of the displacements were obtained by order elevating the functions for the description of the geometry for some of the NURBS patches resulting in the collocation points displayed in \figref{Masjedres}. The material properties of the rock mass, the inclusion and the virgin stresses are given in Table~\ref{tab:practical}. The resulting displaced shape is shown in \figref{Masjedres}.
% \begin{lstlisting}
% Elastic properties of the rock mass:
% E= 10 000 MPa; Poisson's ratio= 0.20

% For the inclusions:
% E=  6 000 MPa; Poisson's ratio= 0.25
% Mohr-Coulomb yield condition, Angle of friction= 30; Cohesion=0.73 MPa; 
% Dilation Angle=0
% \end{lstlisting}
\begin{table}
  \centering
  \begin{tabular}{p{4.5cm}r@{=}l}
    \toprule
    \multicolumn{3}{l}{Rock mass} \\
    Young's modulus   & $E$           & $10000$ MPa    \\
    Poisson's ratio   & $\nu$         & $0.20$        \\ \midrule
    \multicolumn{3}{l}{Inclusion} \\
    Young's modulus   & $E_{i}$       & $6000$ MPa      \\
    Poisson's ratio   & $\nu_{i}$     & $0.25$        \\ \midrule
    \multicolumn{3}{l}{Mohr-Coulomb yield condition} \\
    Angle of friction & $\phi$        & $30^\circ$       \\ 
    Cohesion          & $c$           & $0.73$ MPa \\ \midrule
    Virgin stress field      & $\sigma_{xv}$ &  $-4$ MPa   \\
                      & $\sigma_{yv}$ &  $-8$ MPa   \\ \bottomrule    
  \end{tabular}
  \caption{Practical example: Material parameters and stress field}
  \label{tab:practical}
\end{table}

\section{Summary and Conclusions}

A novel approach to the treatment of piecewise heterogeneous domains with inelastic material behavior has been presented. The approach differs from current methods by the fact that the generation of a cell mesh is replaced by a geometry definition via NURBS and that Kernels of lower singularity are used.

On several test examples, involving finite and infinite domains and different boundary conditions, it is demonstrated that the method works well and leads to accurate results.  It is shown on a practical example that the method is able to handle multiple inelastic inclusions with different material properties.

An important aspect of the implementation is how to deal with nearly singular and singular integration. This involves a large number of Gauss point evaluations and for larger problems needs to be optimized. 

The current software implementation is restricted to tabular inclusions whose geometry can be described with the simple mapping algorithm outlined. However, the methodology for the evaluation of the volume integrals presented here, is not restricted to this description and any geometry definition that allows a mapping to a unit square can be used. It is hoped that the paper will give impetus to much needed further work on this topic.

The paper dealt with plane problems only. An extension of the method and the implementation in 3-D is in progress.

\section{Acknowledgements}
The work was partially supported by the Austrian science fund FWF, under Grant Number P24974-N30: Fast isogeometric boundary element method. This support is gratefully acknowledged.  The basic idea of using body forces instead of initial stresses in the integral equations was arrived at during discussions with John Watson, where the derivation shown in Appendix~\ref{Appendix_A} was also produced.

\appendix
\section{Derivation of integral equation for plasticity}
\label{Appendix_A}

The integral equations for plasticity are derived in tensor notation.
We recall the relationship between stresses and strains
\begin{equation}
  \label{stress}
  \dot{\sigma}_{jk}=C_{jklm}\dot{\epsilon}_{lm}^{e}
\end{equation}
where $\dot{\epsilon}_{lm}^{e}$ is the elastic strain increment and $C_{jklm}$ is the elasticity tensor.

If the elastic limit has been reached, the total strain increment must be split into an elastic and a plastic part
\begin{equation}
  \dot{\epsilon}_{lm}= \frac{1}{2} (\frac{\partial \dot{u}_{l}}{\partial x _{m}} + \frac{\partial \dot{u}_{m}}{\partial x _{l}})= \dot{e}_{lm}^{\epsilon} +  \dot{\epsilon}_{lm}^{p}
\end{equation}
where $\dot{\epsilon}_{lm}^{p}$ denoted the plastic strain increment and $\dot{u}_{l}$ the displacement increment.

Substitution into Equation (\ref{stress}) gives:
\begin{equation}
  \dot{\sigma}_{jk}=C_{jklm}(\dot{\epsilon}_{lm} - \dot{\epsilon}_{lm}^{p})
\end{equation}
or
\begin{equation}
  \label{stress1}
  \dot{\sigma}_{jk}=\dot{\sigma}_{jk}^{e} - \dot{\sigma}_{jk}^{p} 
\end{equation}
where
\begin{equation}
  \dot{\sigma}_{jk}^{e}   =  C_{jklm}\dot{\epsilon}_{lm} 
\end{equation}
and
\begin{equation}
  \dot{\sigma}_{jk}^{p}   =  C_{jklm}\dot{\epsilon}_{lm} ^{p}
\end{equation}
Recall the differential equation of elasticity:
\begin{equation}
  \frac{\partial \dot{\sigma}_{jk}}{\partial x_{k}}  + b_{j}=0
\end{equation}
Substitution of Equation (\ref{stress1}) gives
\begin{equation}
  \frac{\partial \dot{\sigma}_{jk}^{e}}{\partial x_{k}} - \frac{\partial \dot{\sigma}_{jk}^{p}}{\partial x_{k}} + b_{j}=0
\end{equation}
or
\begin{equation}
  \label{diffeb}
  \frac{\partial \dot{\sigma}_{jk}^{e}}{\partial x_{k}} + \dot{b}_{j}^{p} + b_{j}=0
\end{equation}
where the plastic body force $\dot{b}_{j}^{p}= -\frac{\partial \dot{\sigma}_{jk}^{p}}{\partial x_{k}}$ has been introduced.  The stresses are related to the tractions acting on the boundary $S_{0}$ by:
\begin{equation}
  n_{k}\dot{\sigma}_{jk}=\dot{t}_{j}
\end{equation}
or
\begin{equation}
  n_{k}(\dot{\sigma}_{jk}^{e} - \dot{\sigma}_{jk}^{p}) =\dot{t}_{j}
\end{equation}

In terms of $\dot{\sigma}_{jk}^{e}$ we have
\begin{equation}
  n_{k}\dot{\sigma}_{jk}^{e}  =\dot{t}_{j} +\dot{t}_{j}^{p}
\end{equation}
where
\begin{equation}
  \dot{t}_{j}^{p}= n_{k}\dot{\sigma}_{jk}^{p} 
\end{equation}

Applying Betti's theorem we obtain the integral equations (for a full derivation see for example \cite{Banerjee94})
\begin{eqnarray}
  \label{}
  c_{ij} \mal \dot{u}_{i}\left( \pt{y} \right)  = 
  \int_{S} U_{ij}\left( \pt{y},\pt{x} \right) \dot{t}_{j} \left( \pt{x} \right) \dx \boundary  +  \int_{S_{0}} U_{ij} \left( \pt{y},\bar{\pt{x}} \right) \dot{t}_{j}^{p} \left( \bar{\pt{x}} \right) \dx \boundary_{0} \\
  \nonumber
  -   \int_{S}
  T_{ij}\left( \pt{y}, \pt{x}\right) \dot{u}_{j}\left( \pt{x} \right) \dx \boundary +  \int_{V_{0}} U_{ij}\left( \pt{y},\bar{\pt{x}} \right) \dot{b}_{j}^{p}\left( \bar{\pt{x}} \right) \dx \domain_{0} .
\end{eqnarray}

\section{Data for practical example}
\label{Appendix_B}

Definition of excavation boundary:
\begin{lstlisting}
Knot vectors:
0 0 0 1 1 1
0 0 0 1 1 1
0 0 0 1 1 1
0 0 0 1 1 1
0 0 1 1
0 0 1 1
0 0 1 1
0 0 1 1
0 0 1 1
0 0 1 1
0 0 1 1
0 0 1 1

Coefficients:
-15.  30.   0.   1.
-15.  33.1243   0.   0.848 
-12.1919  34.4940   0.   1.

-12.1919  34.494   0.   1.  
-2.5206  38.8109   0.   0.9336 
7.4422  35.9541   0.   1.

7.4422  35.9541   0.   1. 
9.7129  35.3030   0.   0.9962 
11.8360  34.2674   0.   1.0000 

11.8360  34.2674   0.   1. 
15.0000  33.1243   0.   0.8480 
15.0000  30.0000   0.   1. 

15.  30.   0.   1.
15 26 0 1

15 26 0 1
15 21 0 1

15 21 0 1
15 0 0 1

15 0 0 1
-15 0 0 1

-15 0 0 1
-15 10.75 0 1

-15 10.75 0 1
-15 15.75 0 1

-15 15.75 0 1
-15 26 0 1

-15 26 0 1
-15.  30.   0.   1.
\end{lstlisting}

The inclusion were described as follows:
\begin{lstlisting}
Knot vectors:
0 0 1 1
0 0 1 1

0 0 1 1
0 0 1 1

0 0 1 1
0 0 1 1

0 0 1 1
0 0 1 1

Coefficients:
7.4422  35.9541   0   1 
30 53 0 1 
12.1919  34.494   0   1 
31.8 49 0 1 

15 26 0 1
35 37 0 1
15 21 0 1
35 32 0 1

-35 0 0 1
-15 10.75 0 1
-35 5 0 1
-15 15.75 0 1

-35 15 0 1
-15 26 0 1
-35 18 0 1
-15  30   0   1
\end{lstlisting}

\bibliographystyle{elsarticle-num}
\bibliography{ifbbib}

\end{document}